\documentclass[12pt]{amsart}
\pdfoutput=1
\usepackage{geometry}                
\usepackage{graphicx}
\usepackage{amssymb}
\usepackage{epstopdf}
\usepackage[author-year]{amsrefs}  

\newtheorem{theorem}{Theorem}

\newtheorem{lemma}[theorem]{Lemma}
\newtheorem{conjecture}[theorem]{Conjecture}

\def\C{\mathbb C}

\def\Z{\mathbb Z}
\def\T{\mathbb T}

\def\qed{\hfill QED \vskip 3mm}

\def\[{\begin{equation}}
\def\]{\end{equation}}

\newcommand{\comment}[1]{}

\comment{
Files to include with this .tex file
FractalMathieu.bbl  
mathieu50_sym.pdf
mathieu100_central.pdf
mathieu50_warp.pdf
mathieu100_central_top.pdf
mathieu100_central_bot.pdf
mathieu100_side.pdf
mathieu100_cracked.pdf
mathieu100_middle.pdf
Lines_1.pdf
Curves_1.pdf
Lines_1_3.pdf
Curves_1_3.pdf
Curves_1.pdf
Curves_1_3.pdf
Curves_1_5.pdf
Curves_1_7.pdf
Curves_1_9.pdf
Curves_1_2.pdf
Curves_1_4.pdf
Curves_1_6.pdf
Curves_1_8.pdf
Curves_1_3.pdf
Curves_2_3.pdf
CurvePair.pdf
Curves_1_5.pdf
Curves_2_5.pdf
Curves_3_5.pdf
Curves_4_5.pdf
Curves_1_7.pdf
Curves_2_7.pdf
Curves_3_7.pdf
Curves_4_7.pdf
Curves_5_7.pdf
Curves_6_7.pdf
Curves_1_2.pdf
Curves_1_4.pdf
Curves_3_4.pdf
Curves_1_6.pdf
Curves_5_6.pdf
Curves_1_8.pdf
Curves_3_8.pdf
Curves_5_8.pdf
Curves_7_8.pdf
}

\title{Spectra self-similarity for almost Mathieu operators}
\author{Michael P. Lamoureux}
\address[M.~Lamoureux]{Dept.\ Mathematics and Statistics \\ University of Calgary \\ 2500 University Ave NW \\ Calgary AB T2N 1N4 \\ Canada}
\email{mikel@math.ucalgary.ca}
\author{James A. Mingo}
\address[J.~Mingo]{Dept.\ Mathematics and Statistics \\ Queens University \\ Kingston, ON Canada}
\email{mingo@mast.queensu.ca}
\author{Sydney R. Pachmann}
\address[S.~Pachmann]{Faculty of Engineering \\ University of Calgary \\Calgary AB \\ Canada}
\email{sydney.pachmann@gmail.com}
\date{May 7, 2010}                     

\begin{document}

\begin{abstract}
We determine numerically the self-similarity maps for spectra of the almost Mathieu operators, a two-dimensional fractal-like structure known as the Hofstadter butterfly. The similarity maps each have a horizontal component determined by certain algebraic maps, and vertical component determined by a M\"{o}bius transformation, indexed by a semigroup of the matrix group $GL_2(\Z)$. Based on the numerical evidence, we state and prove a continuity result for the similarity maps.  We note a connection between the indexing of the similarity maps and Morita equivalence of rotation algebras $A_\theta$, a continuous field of C*algebras. 
\end{abstract}

\maketitle

\tableofcontents

\section{Introduction}

Looking at the image presented in Figure~1, one immediately notices the striking, repetitive pattern of ``butterfly wings'' that march off towards the vertical horizons at the top and bottom of the image. This rendering of the so-called Hofstadter butterfly, drawn to high-resolution using a combination of numerical algorithms and PostScript graphics programming, reveals some beautiful symmetries of a fundamentally mathematical object. The goal of this paper is to specify exactly the symmetries of this image and prove continuity results of the corresponding similarity maps, motivated by numerical evidence collected in our study of the butterfly.

\begin{figure}[t]
\begin{center}
\includegraphics[width=5in,height=5in]{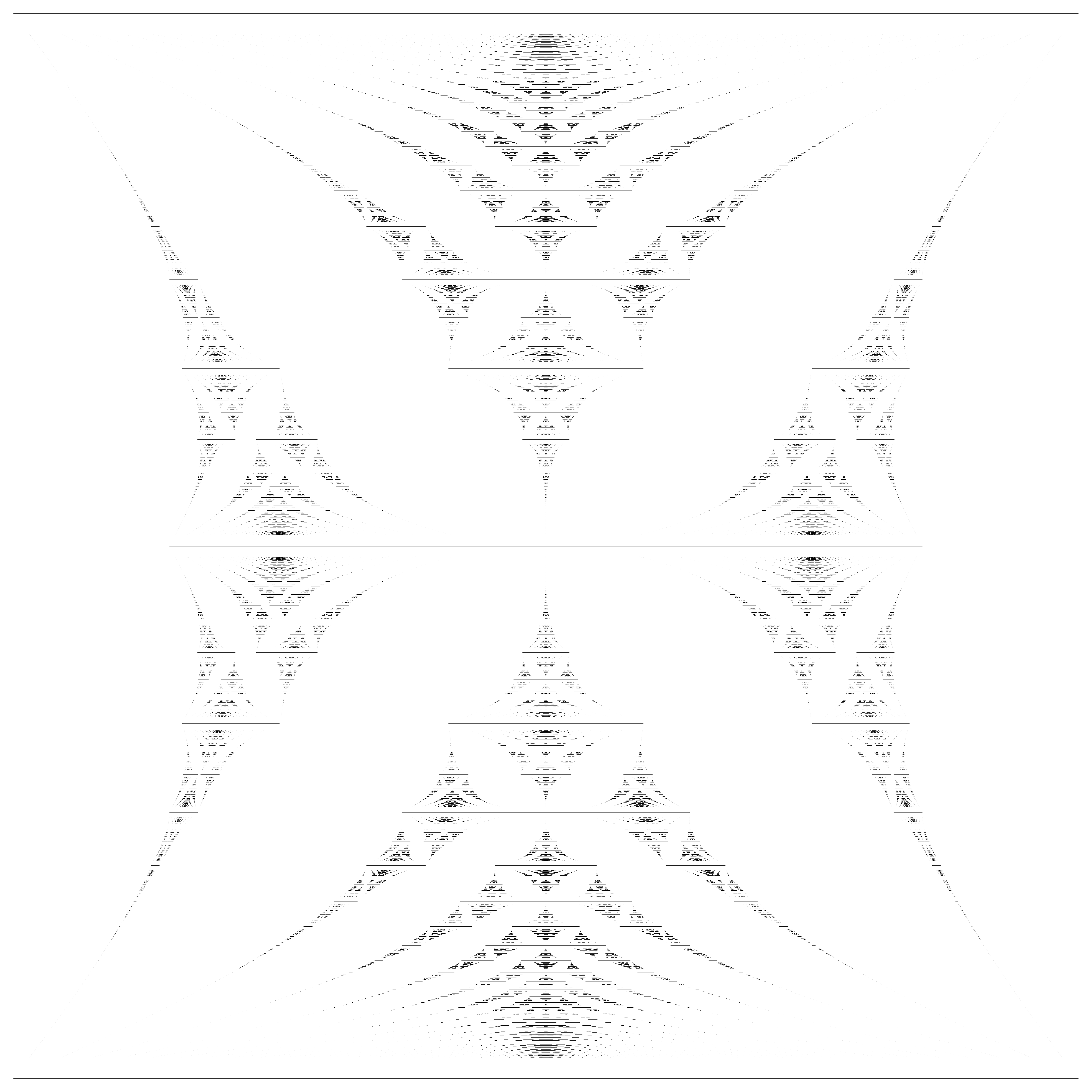}
\caption{The Hofstadter butterfly, a fractal-like structure of line spectra.}
\label{fig:1}
\end{center}
\end{figure}

Although the Hofstadter butterfly looks like a fractal, the image is not created using the usual iterative equations or recursive graphical methods designed to visually render typical fractals~\cite{man90}. Rather, this image is created from an explicit numerical computation of spectral values of a family of linear operators on Hilbert space.

More precisely, this image is constructed as a layered suite of horizontal lines which arise as the spectra of almost Mathieu operators. Each of these operators is succinctly represented as the self-adjoint element 
\[ h_\theta = u + u^* + v + v^* \]
in the C*-algebra $A_\theta$ generated by universal unitaries $u,v$ satisfying the commutation rule
\[ vu = e^{2\pi i \theta} uv. \]
In Figure~1, the vertical axis is spanned by the parameter $\theta$, with $0\leq \theta \leq 1$, while the horizontal axis corresponds to spectral values $x$ in the range $-4\leq x\leq 4$. 

Only rational values of $\theta$ are used in the construction of the image in Figure~1. For this reason, the image is properly called the {\em rational} Hofstadter butterfly. Extensions to irrational values of $\theta$ has been a theme in the long history of study these operators and the physical problems that motivated it, going back at least to a mathematical analysis of Bloch electrons~\cite{brown64}. A selection of relevant studies over the years is given in the references, including~\cites{avila06,bel82,bel90,choi90,connes94,gold09,hof76,kam03,last94,puig03,ypma07}.

Unlike the irrational case as considered in~\cite{arv94}, rational values of $\theta =p/q$ lead to the computation of spectral lines based on finite dimensional eigenvalue calculations using tridiagonal $q\times q$ matrices~\cite{lam97}. In particular, the endpoints of the spectral lines are specified by eigenvalues of these matrices. Numerical algorithms for the tridiagonal eigenvalue computations are rapid and accurate, making calculations of the rational butterfly ideal for numerical experiments. Based on these experiments, we  deduce a systematic catalogue of the symmetries of the butterfly. 

Repetition of a geometric object fading out to infinity is characteristic of hyperbolic geometry~\cite{cox42}, hence it is perhaps no surprise that the matrix group $GL_2(\Z)$ arises in the observed symmetries. Sections~2 through 5 reveal a semigroup of linear fractional transformations represented by elements of $GL_2(\Z)$ which act on the vertical parameter $\theta$ to generate the family of similarity maps.

Equally important are the algebraic curves specifying the horizontal component of the similarities, acting on the horizontal parameter $x$. Sections~6 through 9 present numerical evidence of the continuous maps taking one level of horizontal spectra to another, using an indexing of intervals and a correspondence of image points under polynomial maps. These are precisely the characteristic polynomials of the aforementioned finite dimensional $q\times q$ matrices giving the spectra.

Uniting the vertical and horizontal components leads to an identification of the general form of the similarity maps, as presented in Section~10. The similarity extends to gap-labelling, a specific method for indexing the butterfly wings based on Chern characters, which is discussed in Section~11. Proof of continuity of the similarity maps is given in Section~12 and generators for the semigroup of similarities is presented in Section~13. 

``Xenocides, who is ugly, makes ugly poetry," said Aristophanes. With the Hofstadter butterfly, we see pretty mathematics making pretty pictures.

\section{Numerical evidence: Similarity maps}

Mirror image symmetry in the horizontal direction suggests an obvious self-similarity of the image in Figure~1  given by a reflection about the line $x=0$, namely:
\[ (x,\theta) \mapsto (-x, \theta). \]
This map is continuous, and does map the butterfly to itself, as the spectrum of each $h_\theta$ is symmetric. 

In the vertical direction, another mirror image symmetry in the figure suggests a second self-similarity map given by reflection about the line $\theta = 1/2$, namely:
\[ (x,\theta) \mapsto (x, 1-\theta). \]
Again, this is a continuous map, and again maps spectra properly since the algebra $A_\theta$ is isomorphic to $A_{1-\theta}$, with operator $h_\theta$ mapping onto $h_{1-\theta}$ under the isomorphism.

Now, a more interesting symmetry is observed mapping the large central butterfly onto the next largest butterfly in the bottom half of the image in Figure~1. Zooming in on this butterfly, as shown in Figure~2, one sees the top of the butterfly at $\theta = 1/3$, the bottom at $\theta = 0$, and the centre of the butterfly at $\theta = 1/4$. This suggest we must find a self-similarity map that, on vertical parameter values $\theta$, will map $\theta \mapsto \theta'$ as
\[ 0 \mapsto 0, \qquad 1/2 \mapsto 1/4, \qquad 1 \mapsto 1/3. \]
No linear map will do, but a linear fractional transformation\footnote{i.e. a M\"{o}bius transformation} does work, using the map
\[ \theta \mapsto \theta' = \frac{\theta}{2\theta + 1}. \]
It is reassuring to notice that other obvious lines in the full image map correctly under this map. For instance, the line at $\theta = 1/3$ in the full butterfly in Figure~1 should map to the top of the second largest butterfly in Figure~2, which is at $\theta' = 1/5.$ And indeed, the LFT here does just that. Similarly we can check for that the line at $1/4$ maps to $1/6$, the line at $1/5$ maps to $1/7$, and so on. Checking numerically many of these lines assures us that the linear fractional transformation is the proper choice.

\begin{figure}[t]
\begin{center}
\includegraphics[width=5in]{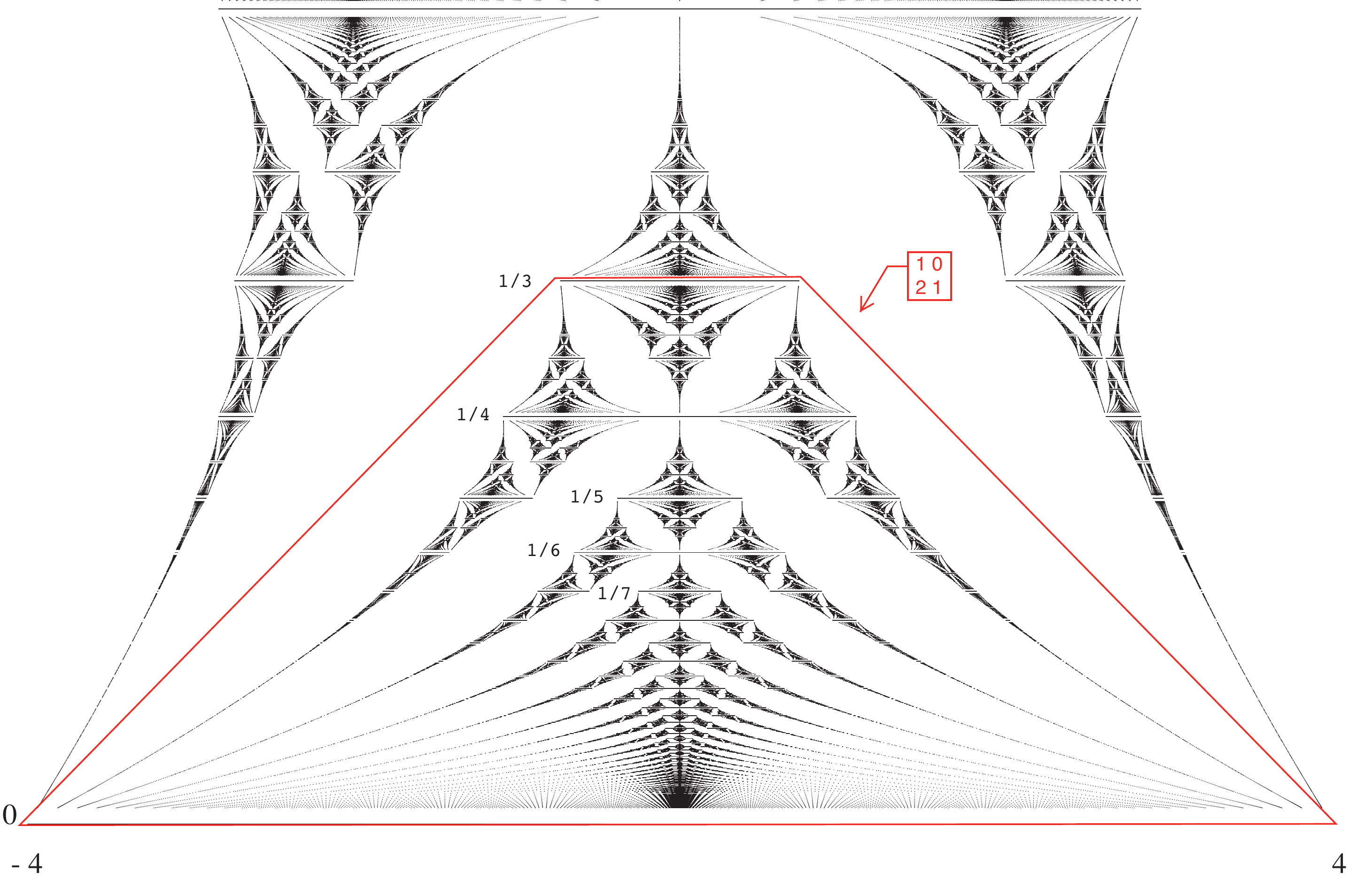}
\caption{Similarity to the bottom third of the butterfly.}
\label{fig:2}
\end{center}
\end{figure}

Going to a numerical experiment based on these observations, we construct a map from the full image in Figure~1 to subset like Figure~2, using the LFT in the vertical component, and contracting linearly in the horizontal component. Specifically, we map according to the rule
\[ 
(x,\theta) \mapsto (x',\theta') = \left((1-.82*\theta)x ,\frac{\theta}{2\theta + 1} \right). 
\]
The coefficient $0.82$ was chosen to get the correct width for the top horizontal line. The result is the image shown in Figure~3, which is very much like the lower butterfly image in Figure~2. 

Observing the slight differences between the image in Figure~1 and the result of the numerical experiment shown in Figure~3, (the difference in the curvature of the wingtips for instance), we can conclude the whatever the similarity map is, it is only approximately linear. Nevertheless, the basic structure and position of the tiny butterflies are preserved. This serves as numerical evidence that the vertical maps are indeed given by linear fractional transformations. The horizontal component is more complicated, and is discussed in Sections~7, 8, 9.

\begin{figure}[t]
\begin{center}
\includegraphics[width=5in]{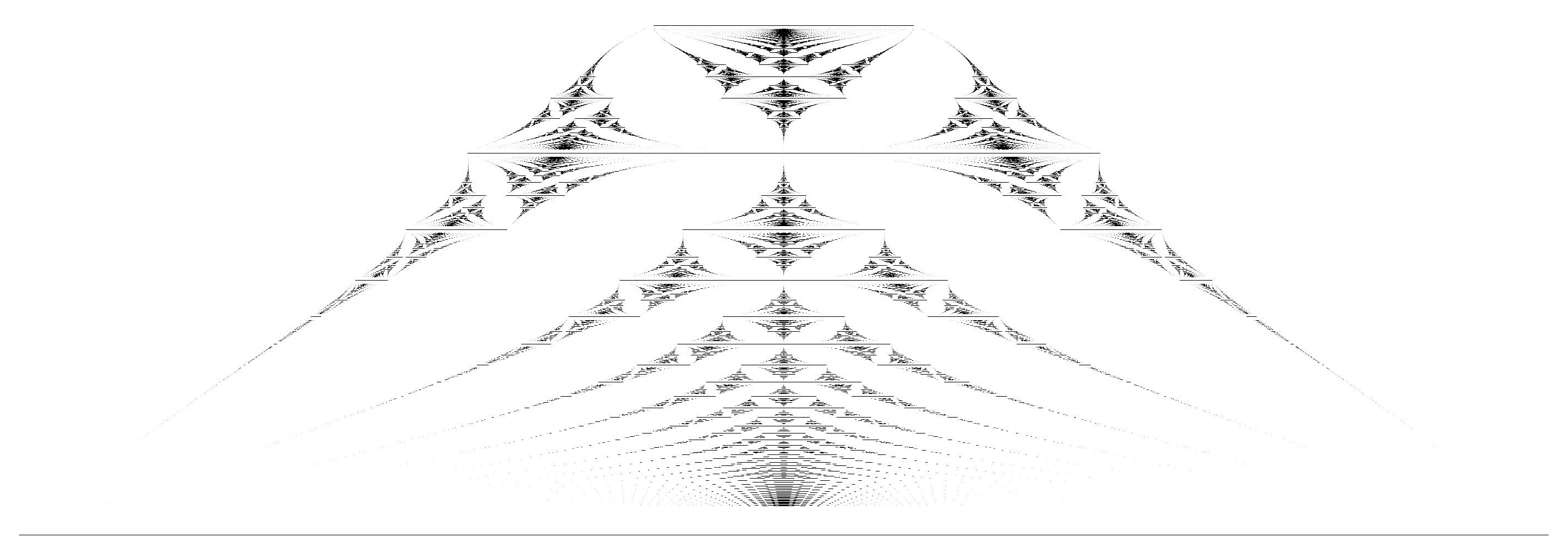}
\caption{A rendering of an approximate similarity to the bottom third, using a linear fractional transformation vertically, linear scaling horizontally.}
\label{fig:3}
\end{center}
\end{figure}

\section{Numerical evidence: Vertical similarities from $GL_2(\Z)$}

Perhaps it is suggestive that two of the self-similarity maps identified so far (vertical flip, and the bottom third map) have vertical components give by M\"{o}bius transformations,
\[ \theta \mapsto \frac{1-\theta }{1} \mbox{ and } 
\theta \mapsto \frac{\theta}{2\theta + 1}. \]
Both transformations are specified by a $2\times 2$ matrix in $GL_2(\Z)$, in the form
\[ 
M =
\left[
\begin{array}{rr}
a & b \\
c & d
\end{array}
\right],
\]
yielding a corresponding linear fractional transformation 
\[ \theta \mapsto \theta' = \frac{a\theta +b}{c\theta + d}. \]
It will be convenient to label these maps by the corresponding matrix, and we note that composition of the LFTs corresponds to matrix multiplication in $GL_2(\Z)$.

\begin{figure}[t]
\begin{center}
\includegraphics[width=5in]{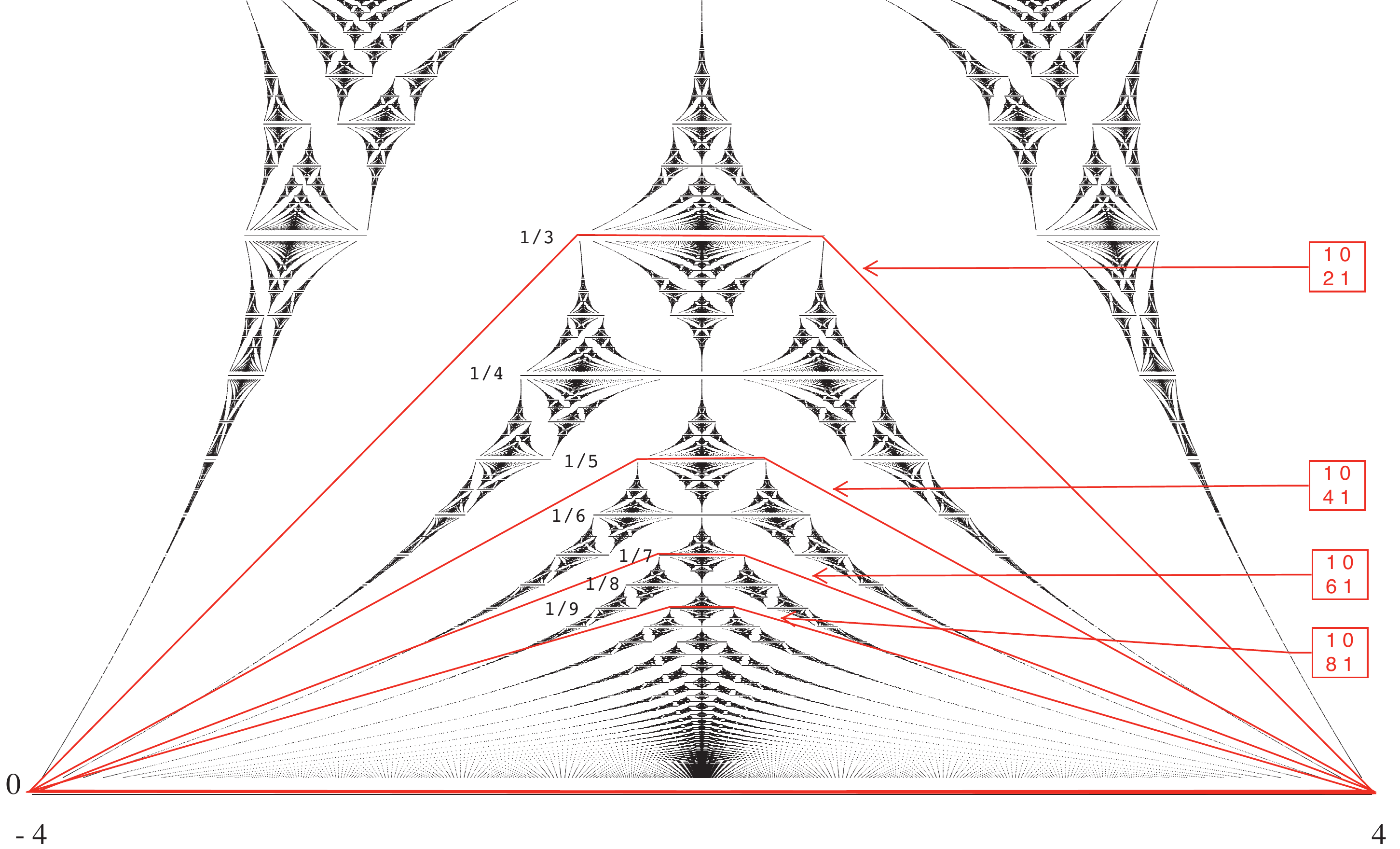}
\caption{A sequence of similarities in the lower central core of the butterfly, and matrices in $SL_2(\Z)$ that implement the vertical component of the similarity.}
\label{fig:4}
\end{center}
\end{figure}

A leap of faith suggests looking for symmetries in the butterfly indexed by elements of $GL_2(\Z)$. 
Examining Figure~4, we see the bottom of the central core of butterflies, and we can compute a list of matrices that implement the self-similarity map on the vertical component.
For instance, we can map the whole butterfly onto the central bottom butterfly that tops out at label $\theta = 1/5$. For this, we need an LFT that maps $0\mapsto 0, 1/2 \mapsto 1/6, 1\mapsto 1/5$. Solving for coefficients $a,b,c,d$ in the LFT, we find the matrix
\[ 
\left[
\begin{array}{rr}
1 & 0 \\
4 & 1
\end{array}
\right]
\]
will work to implement an LFT that maps the whole figure into this central bottom butterfly.

Continuing down the central core, we see an infinite sequence of butterflies extending to the ``horizon'' at $\theta = 0$. Some simple calculations analogous to those above produce a sequence of  similarity maps indexed by matrices of the form
\[ 
\left[
\begin{array}{ll}
1 & 0 \\
2n& 1
\end{array}
\right] =
\left[
\begin{array}{ll}
1 & 0 \\
1& 1
\end{array}
\right]^{2n} = A^{2n}, \mbox{ for $n \geq 1$. }
\]
We will see later that this matrix $A$ is a key generator of the similarities. 

\begin{figure}[t]
\begin{center}
\includegraphics[width=5in]{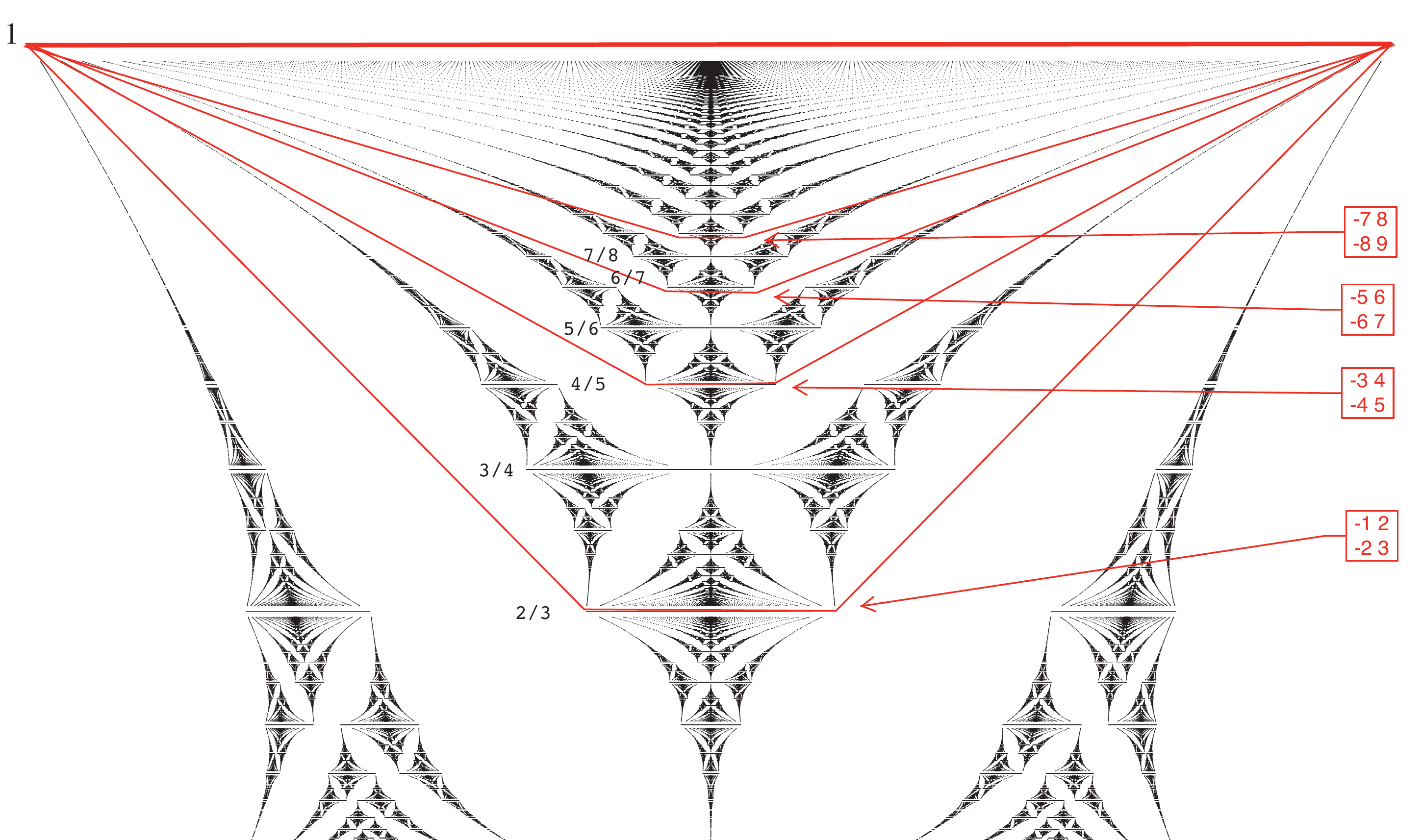}
\caption{A sequence of similarities in the upper central core of the butterfly, and matrices in $GL_2(\Z)$ that implement the vertical component of the similarity.}
\label{fig:5}
\end{center}
\end{figure}

Hopping to the top of butterfly, in Figure~5 we see a central core of butterflies extending to the upper ``horizon'' at level $\theta=1$. Here, the self-similarity maps are obtained using matrices in $GL_2(\Z)$ of the form
\[ 
\left[
\begin{array}{ll}
1-2n & 2n \\
-2n& 2n+1
\end{array}
\right], \mbox{ for $n \geq 1$. }
\]

Matrices from symmetries on top (Eqn~13)  are related to matrices for the symmetries on the bottom (Eqn~12), via the conjugation
\[ 
\left[
\begin{array}{ll}
1-2n & 1 \\
-2n & 2n+1
\end{array}
\right] = 
\left[
\begin{array}{rr}
-1 & 1 \\
 0 & 1
\end{array}
\right]
\left[
\begin{array}{ll}
1 & 0 \\
2n & 1
\end{array}
\right] 
\left[
\begin{array}{rr}
-1 & 1 \\
 0 & 1
\end{array}
\right].
\]
This is just conjugation with the flip symmetry map of the butterfly, as the matrix
\[ B = 
\left[
\begin{array}{rr}
-1 & 1 \\
 0 & 1
\end{array}
\right]
\]
induces the linear fractional transformation mapping $\theta$ to $1-\theta$, turning the butterfly upside down.

\begin{figure}[t]
\begin{center}
\includegraphics[width=3in]{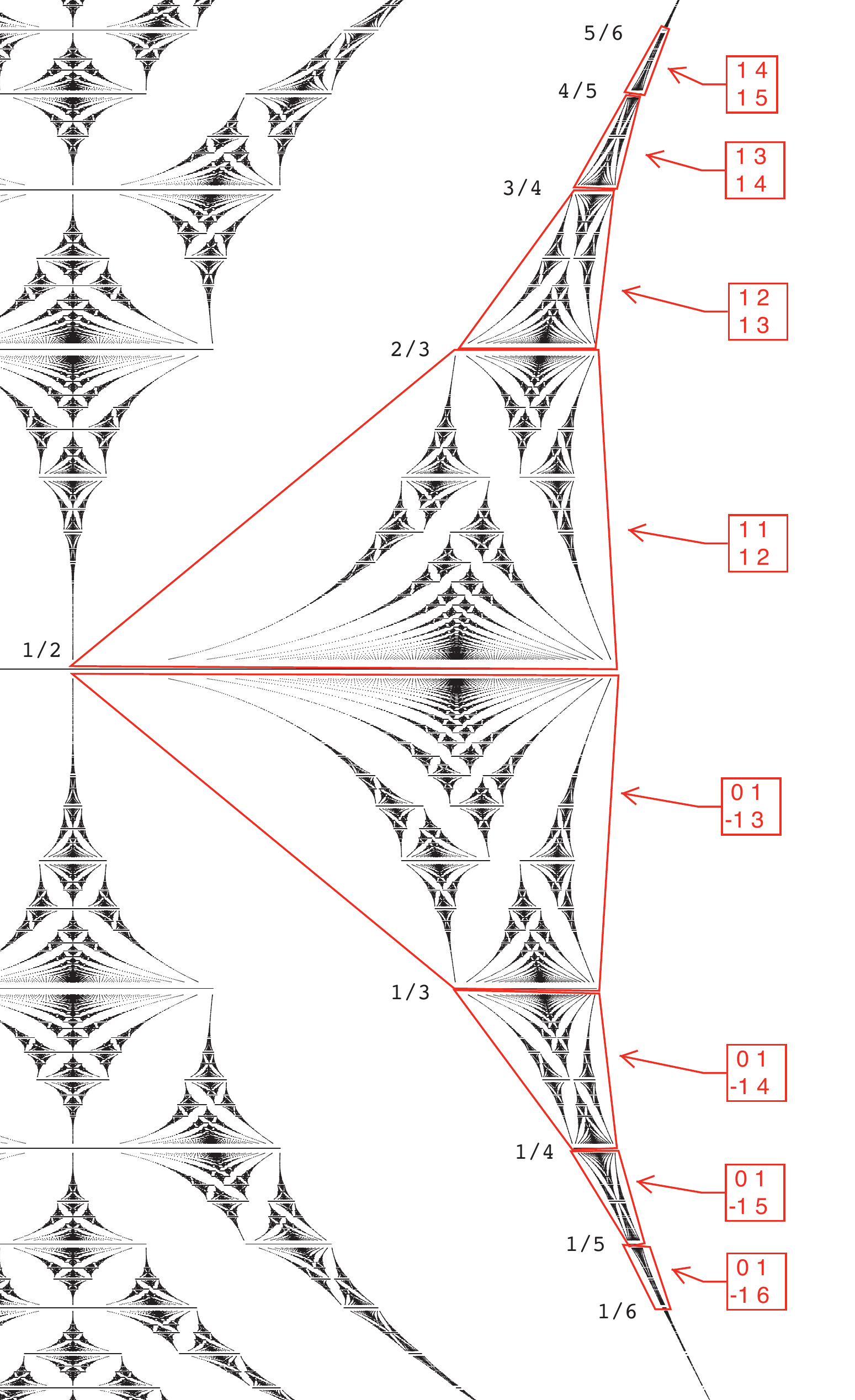}
\caption{A sequence of similarities on the edge of the butterfly, and matrices in $GL_2(\Z)$ that implement the similarity.}
\label{fig:6}
\end{center}
\end{figure}

Attending to some of the butterflies on the side of the image, we find the relevant linear fractional transformations are represented by matrices of the form
\[
\left[
\begin{array}{rr}
0 & 1 \\
-1 & n
\end{array}
\right] , \mbox{ for $n \leq 3$},
\]
which generate the symmetries on the bottom half of the image shown in Figure~6. On the top of the image, we obtain transformations indexed by matrices of the form
\[
\left[
\begin{array}{ll}
1 & n-2 \\
1 & n-1
\end{array}
\right] , \mbox{ for $n \geq 3$}.
\]
Again, we observe there is a relation between these two forms of transformation, through conjugation with the flip, as we have
\[ 
\left[
\begin{array}{ll}
1 & n-2 \\
1 & n-1
\end{array}
\right] = 
\left[
\begin{array}{rr}
-1 & 1 \\
 0 & 1
\end{array}
\right]
\left[
\begin{array}{rr}
0 & 1 \\
-1 & n
\end{array}
\right] 
\left[
\begin{array}{rr}
-1 & 1 \\
 0 & 1
\end{array}
\right].
\]

Note, however, that these self-similarity maps have a certain discontinuity. In particular, the centre of the large butterfly gets mapped to a ``broken'' butterfly on these side maps. 
This tells us the self-symmetry maps need not be continuous. In particular there may be a discontinuity in the horizontal direction as we cross the $x=0$ spectral value.

Nevertheless, this is a mild discontinuity and is easy to understand. In the construction of the Hofstadter butterfly,  for $\theta = p/q$ with $q$ even, the spectrum consists of $q$ intervals on the real line, two of which touch at the point zero.~\cite{choi90} It is these ``touching intervals'' that are getting broken apart in the above similarity maps. So although there is an apparent discontinuity, it is perhaps better to describe it as the breaking apart of two touching spectral lines. The similarity map should take the single, common endpoint, and map it to two distinct endpoints of two non-overlapping intervals. The details of this double-valued map will be discussed in Section~10.

\section{Numerical evidence: A symmetry with break}

As noted in the previous section, under the similarity map sometimes the butterfly breaks. This suggests we can look for more similarities if we are a bit more open to what a broken butterfly looks like. 

\begin{figure}[t]
\begin{center}
\includegraphics[width=5in]{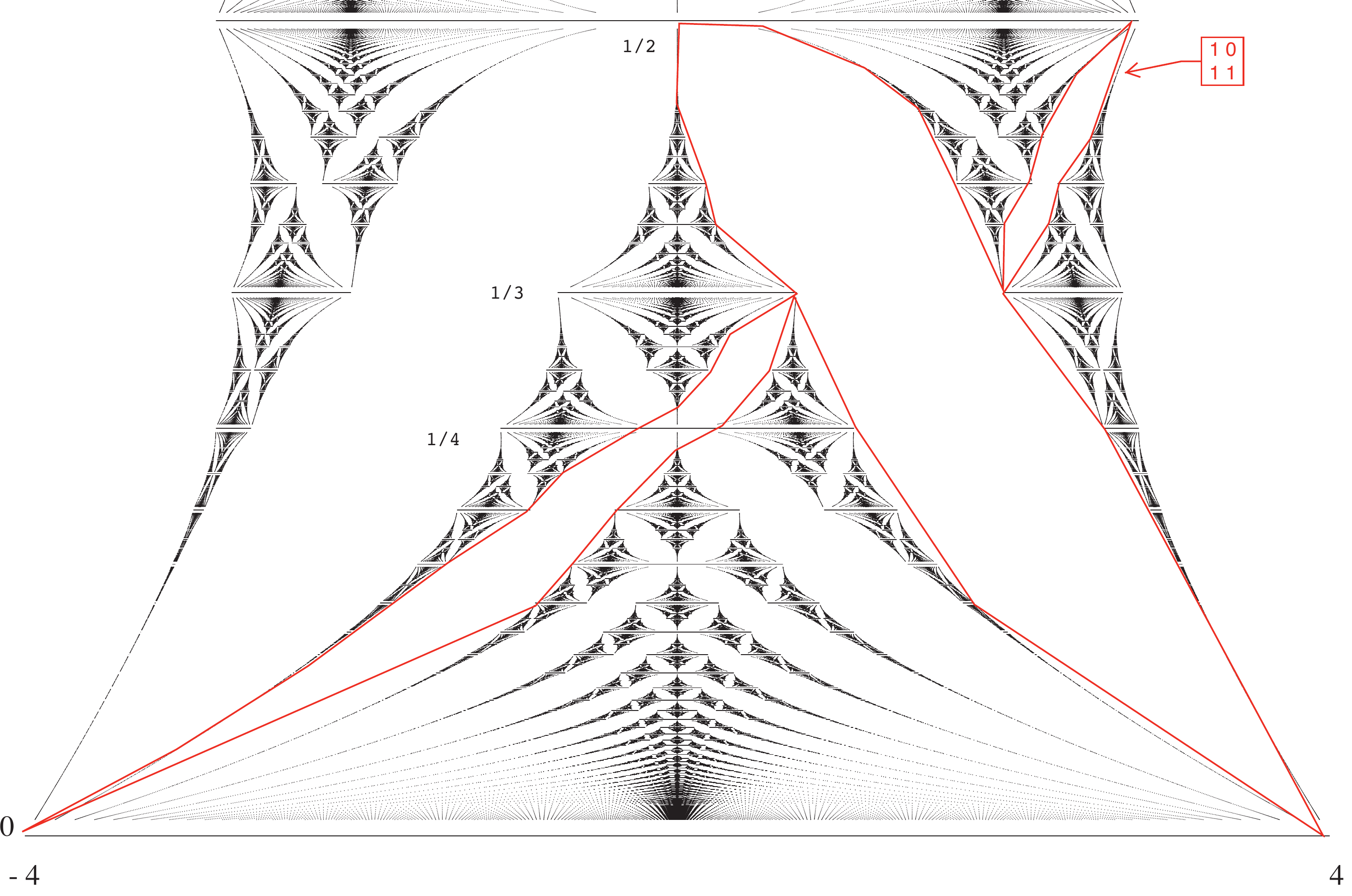}
\caption{An unusual, cracked similarity. Note the basic butterfly structure spanning the bottom half of the full image.}
\label{fig:7}
\end{center}
\end{figure}

In Figure~7, we have an example of a really broken butterfly. The butterfly fits into the region
\[ 0\leq \theta \leq 1/2, \]
and the similarity map is given by the linear fractional transformation with matrix
\[
\left[
\begin{array}{rr}
1 & 0 \\
1 & 1
\end{array}
\right].
\]
This break (at $\theta = 1/3$) looks pretty bad, but in fact we know from gap labelling theorems that when one plots the gaps an inverse slope 2, there is a discontinuity exactly at $\theta = 1/3$~\cite{lam10}. So in fact this is just confirming the fact that the butterflies are coming from gap labelling. 

By including this symmetry, we include the matrix
\[
A = \left[
\begin{array}{rr}
1 & 0 \\
1 & 1
\end{array}
\right]
\]
in our collection of matrices in $GL_2(\Z)$ generating similarity maps. 

We observe that all the LFTs seen so far map the $\theta$ interval $[0,1]$ into itself, and are indexed by certain elements of of $GL_2(\Z)$. This evidence suggests a certain semigroup in $GL_2(\Z)$ specifies the possible self-similarity maps. A precise statement identifying the semigroup is given in the next section.

\section{Generating the $GL_2(\Z)$ symmetries}

The following technical result states that the set of linear fractional transformations which map the interval $[0,1]$ into itself forms a semigroup, indexed by a 2-generated semigroup of  matrices in $GL_2(\Z)$.

\begin{theorem}
Let $G$ be the semigroup of linear fractional transformations of the form
\[ \theta \mapsto \theta' = \frac{a\theta + b}{c\theta +d}, \qquad  \mbox{ for some }
\left[
\begin{array}{rr}
a & b \\
c & d
\end{array}
\right] \in GL_2(\Z)/\{\pm I\},
\]
which map the interval $[0,1]$ into itself. Then $G$ is generated by the two maps $\phi(\theta) = \theta/(\theta+1)$ and $\psi(\theta) = 1-\theta$, represented by matrices
\[ A = 
\left[
\begin{array}{rr}
1 & 0 \\
1 & 1
\end{array}
\right], \mbox{ and }
B = \left[
\begin{array}{rr}
-1 & 1 \\
0 & 1
\end{array}
\right].
\] 
\end{theorem} 

The proof of this technical result is given in Appendix~1, and is based on the Euclidean algoritm. We make a few observations about this theorem.

The LFTs do not notice a change of sign in the representing matrix, so the theorem refers to the quotient group $GL_2(\Z)/\{\pm I\}$.

These two generators correspond to two symmetries for the Hofstadter butterfly. The first generates the map that takes the butterfly to its lower half, mentioned in Section~4. The second is the vertical flip symmetry, mentioned in Section~2. 
By including the vertical flip in our symmetries,  we are able to generate all (vertical components) of the symmetries using only two generating maps. 

The fact that the Euclidean algorithm is involved suggests that determining the details of a symmetry map will be rather involved -- the factorization into generators is not simple.

It is worth mentioning a deep result from C*-algebras, that states two rotation algebras  $A_\theta, A_{\theta'}$ are Morita equivalent if and only if there is a linear fractional transformation mapping
\[ \theta \mapsto \theta' = \frac{a\theta +b}{c\theta + d}, \]
for some matrix 
$\left[
\begin{array}{rr}
a & b \\
c & d
\end{array}
\right] 
$
in $GL_2(\Z)$~\cite{effros80}. We have no idea what connection this might have with the self-similarity maps we have above, which act on spectra of operators in the algebras $A_\theta, A_{\theta'}$, not on the algebras themselves. 

\vskip 5mm

\section{Horizontal similarity: Interval maps}

Considering now the horizontal component of the similarity maps, a close examination of the apparent symmetries suggests the intervals $I_1, I_2, \ldots, I_q$ in the spectrum at level $\theta = p/q$ are mapped onto a subcollection of the intervals $I'_1,  I'_2, \ldots I'_{q'}$ in the spectrum at level $\theta' = p'/q'$.  We note here which subcollections appear.

\begin{figure}[t]
\begin{center}
\includegraphics[width=5in]{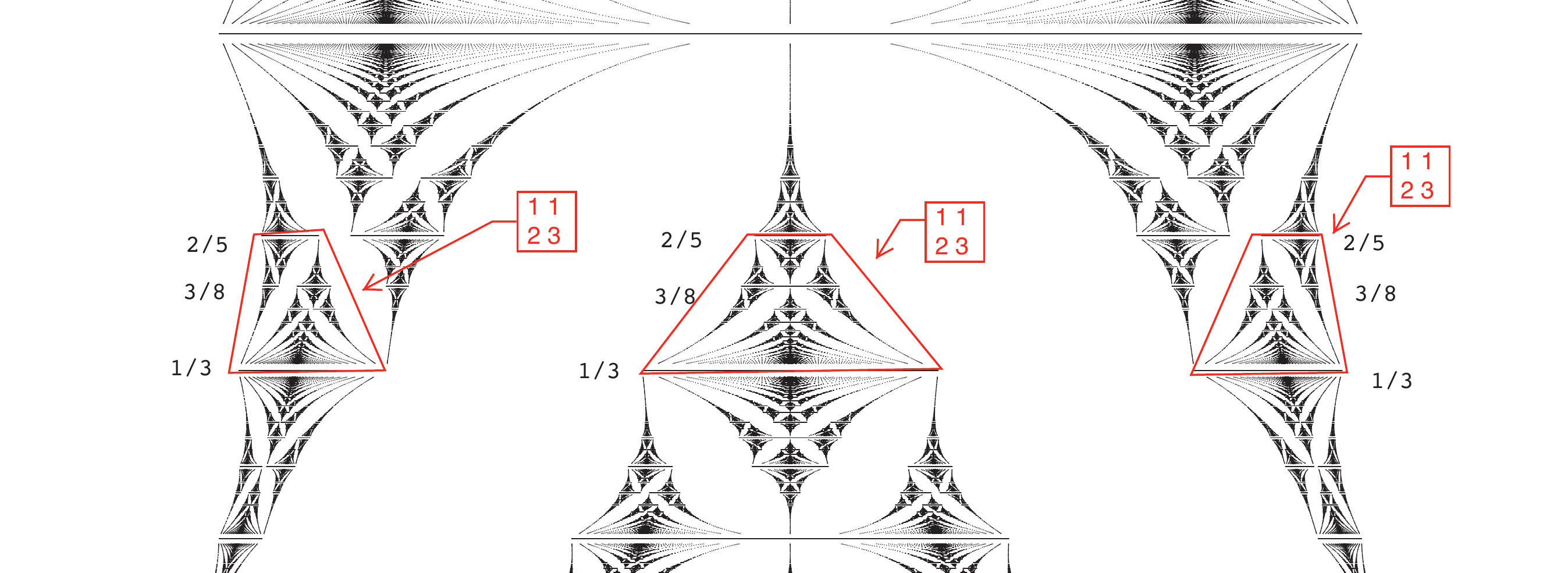}
\caption{Multiple similarities, indexed by the same matrix in $SL_2(\Z)$.}
\label{fig:8}
\end{center}
\end{figure}

Examining Figure~8, notice for the single vertical map with linear fractional transformation given by matrix
\[ 
\left[
\begin{array}{rr}
1 & 1 \\
2 & 3
\end{array}
\right],
\]
there are in fact three choices of horizontal maps. Start first at initial $\theta = 1$, which maps to $\theta' = 2/5$ under this LFT. The spectrum at level $\theta=1$ has one interval, while at  level $\theta' = 2/5$, there are five intervals that we can map to. Only three are obtained. In the left part of the image in Figure~8, the initial interval $I_1$ maps onto interval $I'_1$; the middle image maps onto interval $I'_3$, and the right image maps onto interval $I'_5$. Writing $p'/q' = 2/5$, we see that the three interval maps are given as
\[ 
I_1 \mapsto I'_1 = I'_{0p' + 1}, \qquad 
I_1 \mapsto I'_3 = I'_{1p' + 1}, \qquad 
I_1 \mapsto I'_5 = I'_{2p' + 1}. \]

Checking now at the level $\theta=1/2$ mapping to $\theta'=3/8$, the two intervals $I_1, I_2$ at $\theta=1/2$ map to two interval $I'_1, I'_2$ on the left; to the two interval $I'_4,I'_5$ in the middle; and $I'_7,I'_8$ on the right. In this case, we have $p'/q'=3/8$ and we can summarize the left, middle and right maps as
\[ 
I_k \mapsto  I'_{0p' + k}, \qquad 
I_k \mapsto  I'_{1p' + k}, \qquad 
I_k \mapsto  I'_{2p' + k}, \]
for $k=1,2$. Similarly, at level $\theta = 0$, the single interval $I_1$ maps as
\[ 
I_1 \mapsto I'_1 = I'_{0p' + 1}, \qquad 
I_1 \mapsto I'_3 = I'_{1p' + 1}, \qquad 
I_1 \mapsto I'_5 = I'_{2p' + 1}, \]
where $p'/q' = 1/3$. The form of the interval map is very consistent.

We go back and check the other similarity maps. In Figure~4, the interval maps are apparently of the form
\[ I_k \mapsto  I'_{n\cdot p + k} = I'_{n\cdot p'  + k}, \]
for the vertical similarity given by 
\[ 
\left[
\begin{array}{rr}
1 & 0 \\
2n & 1
\end{array}
\right].
\]
We can see this by observing the $q$ intervals at level $\theta = p/q$ map to the central $q$ intervals of the set of $q' = 2np + q$ intervals at level $\theta' = p'/q$.

In Figure~5, the interval maps are of the form
\[ I_k \mapsto  I'_{n\cdot (q-p) + k} = I'_{n\cdot (q'-p')  + k}, \]
for the vertical similarity given by 
\[ 
\left[
\begin{array}{rr}
1-2n & 2n \\
-2n & 2n+1
\end{array}
\right].
\]
Again, the $q$ intervals at level $\theta = p/q$ map to the central $q$ intervals of the set of $q' = 2n(q-p) + q$ intervals at level $\theta' = p'/q$.

In Figure~6, a similarity in the top half, say with matrix
\[ 
\left[
\begin{array}{rr}
1 & 1 \\
1 & 2
\end{array}
\right],
\]
gives interval maps of the form
\[ I_k \mapsto  I'_{r\cdot p' + k}, \]
with $r=1$ on the right side of the butterfly (as illustrated), and $r=0$ on the left side of the butterfly (not shown in Figure~6). However, a similarity in the bottom half of Figure~6, such as with matrix
\[ 
\left[
\begin{array}{rr}
0 & 1 \\
-1 & 3
\end{array}
\right],
\]
gives interval maps of the form
\[ I_k \mapsto  I'_{r\cdot (q'-p') + k}, \]
with $r=1$ on the right map, $r=0$ on the left map.

Generally speaking, it appears there are two types of interval maps, those of the form
\[ I_k \mapsto  I'_{r\cdot p' + k}, \]
and those of the form
\[ I_k \mapsto  I'_{r\cdot (q'-p') + k}, \]
for some choices of non-negative integer $r$.

In the next three sections we determine how the horizontal component of the similarity map behaves within each interval $I_k$.

\section{Horizontal similarity: Cubic case $1\mapsto 1/3$}
The vertical component of the self-similarity map is given by linear fractional transformations, as described in Sections~2 through 4. The horizontal component maps intervals to intervals in the spectra, and appears to be nearly linear, as we saw in the construction of Figure~3. We investigate the details of this horizontal map. 

Again referring to the butterfly list in Figure~1, we see that the spectral line at height $\theta = 1$ is approached by a dense cluster of spectral points a heights $\theta <1$. Under the LFT mapping $\theta \mapsto \theta' = \theta/(2\theta +1)$, the point $\theta=1$ maps to $\theta' = 1/3$, and this dense cluster of spectral points near 1 maps to a dense cluster near $1/3$. 

We would expect by continuity arguments that the horizontal map at $\theta=1$ should be well-approximated by following how the nearby cluster points map. Setting rational $p/q = \theta <1$, we have $p \approx q$, and under the LFT we have $p'/q' = p/(2*p+q)$. There are $q$ points near line 1, and they map to the middle third of the $2p+q \approx 3q$ points under the LFT.

So we do a numerical experiment: Map these cluster points near 1 to corresponding points near $1/3$, and plot the result, to obtain Figure~9. We note the curve is almost linear, but not quite. Performing a best fit polynomial approximation\footnote{In fact, in our first numerical experiments, we inadvertently computed the inverse of a polynomial, with expansion $y = x/6 + x^3/1296 + x^5/93312 + \cdots$. The second author recognized this as the inverse of a cubic!}, we find the map is in fact given by
\[ x = 6y - y^3. \]

We observe that these polynomial maps are known from earlier work~\cite{choi90} as characteristic polynomials $P_\theta$ for matrices used to compute the spectra of $q\times q$ matrices representing the operators $h_\theta$ in the rotation algebra $A_\theta$, with $\theta = p/q$. Thus, the map we observe is simply
\[ P_1(x) = -P_{1/3}(y), \]
where by $P_{p/q}(x)$ we mean the q-th order characteristic polynomial of the $q \times q$ matrix
\[ H = U + U^* + V + V^*,\]
with $U$ the cyclic permutation matrix and $V$ a diagonal matrix with consecutive powers of $e^{2\pi i p/q}$ on the diagonal. Details are in \cites{choi90, lam07}.

To see the full algebraic curve, we plot the full curve in Figure~10. Observe the central portion of the curve corresponds to the nearly linear section we saw in Figure~9.

\begin{figure}[t]
\begin{center}
\includegraphics[width=3in]{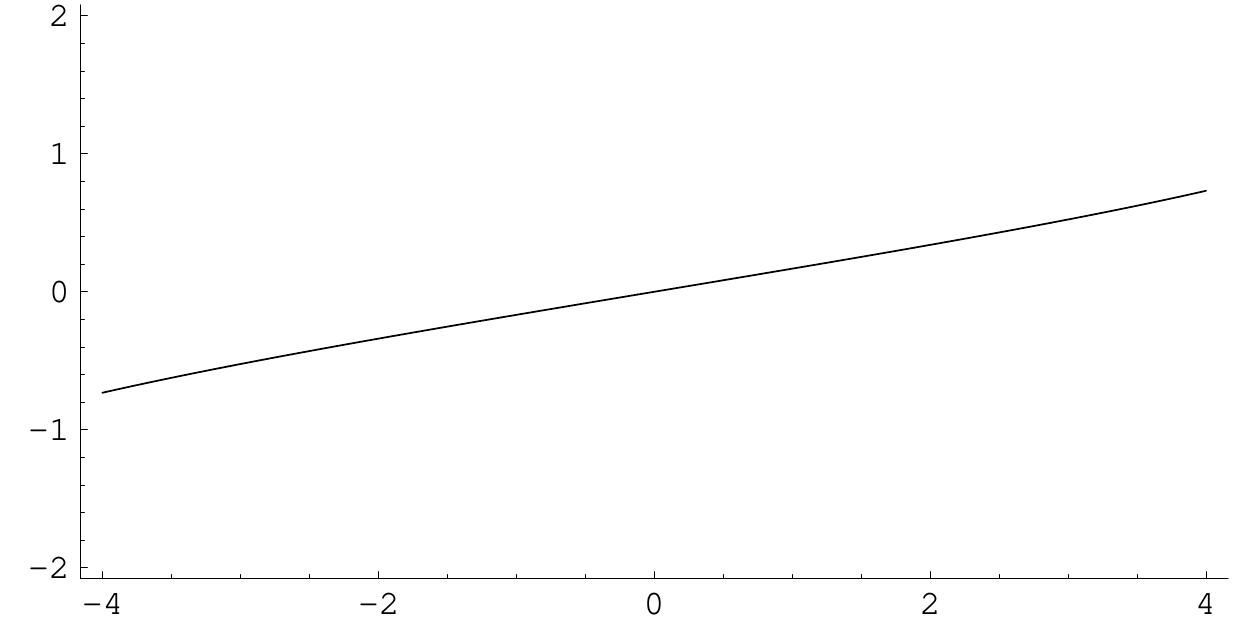}
\caption{Map from $\theta=1$ line to 1/3 lines (middle interval).}
\label{fig:9}
\end{center}
\end{figure}

\begin{figure}[t]
\begin{center}
\includegraphics[width=3in]{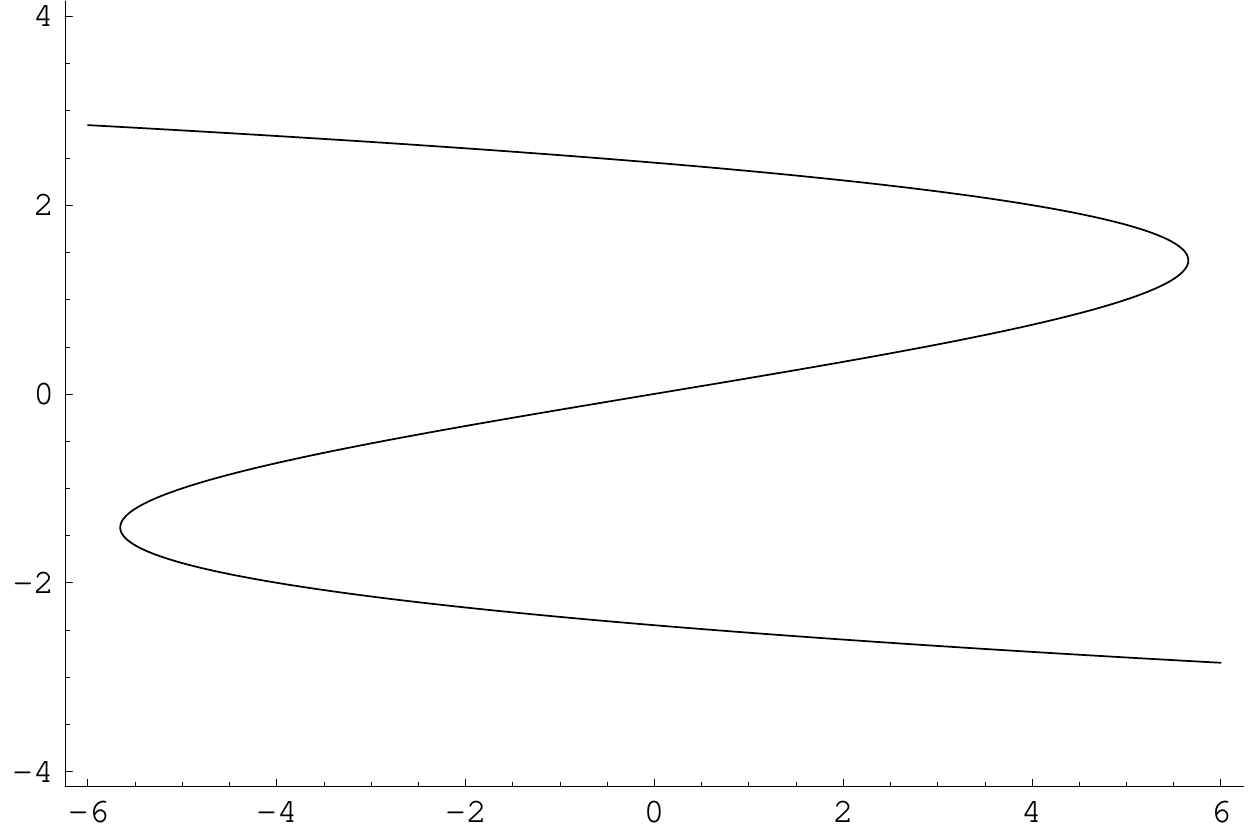}
\caption{Algebraic curve $P_1(x) = -P_{1/3}(y)$ defines the 1 to 1/3 map.}
\label{fig:10}
\end{center}
\end{figure}

\section{Horizontal similarity: Quintic case $1/3 \mapsto 1/5$}

We repeat the numerical calculation, using the same central symmetry with LFT given by the matrix
\[ 
\left[
\begin{array}{rr}
1 & 0 \\
2 & 1
\end{array}
\right].
\]
But now we look how the three spectral lines at $\theta=1/3$ maps to the three central spectral lines  at level $\theta' = 1/5$. A numerical experiment as in the last section produces the plot in Figure~11. Again, the three segments look nearly linear, but not quite. An inspired guess suggest these segments come from a portion of the algebraic curve
\[  P_{1/3}(x) = -P_{1/5}(y), \]
and a comparison of the two plots shows this is indeed the case.

\begin{figure}[t]
\begin{center}
\includegraphics[width=3in]{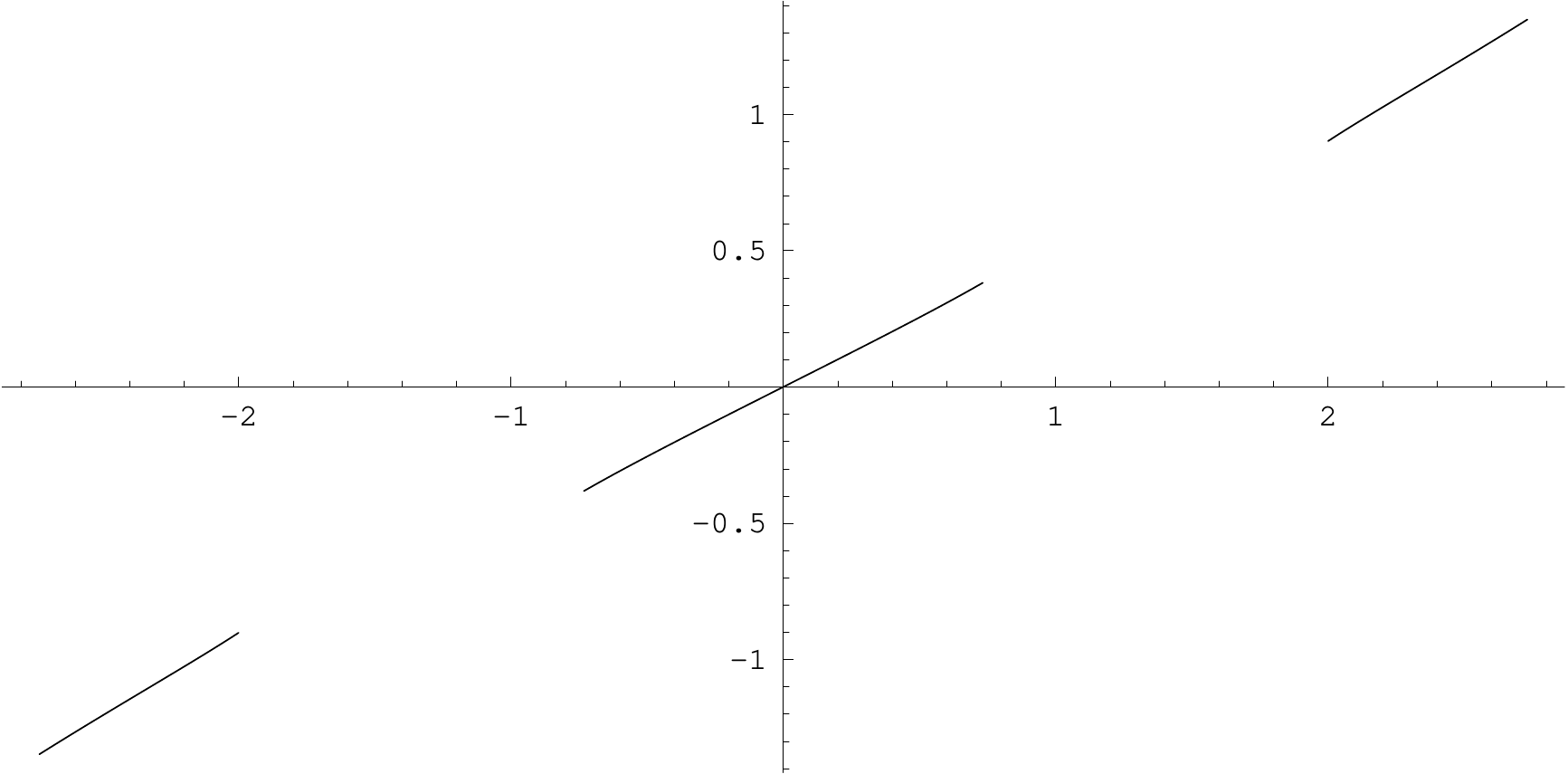}
\caption{Map from 1/3 lines to 1/5 lines (3 middle intervals).}
\label{fig:11}
\end{center}
\end{figure}

\begin{figure}[t]
\begin{center}
\includegraphics[width=3in]{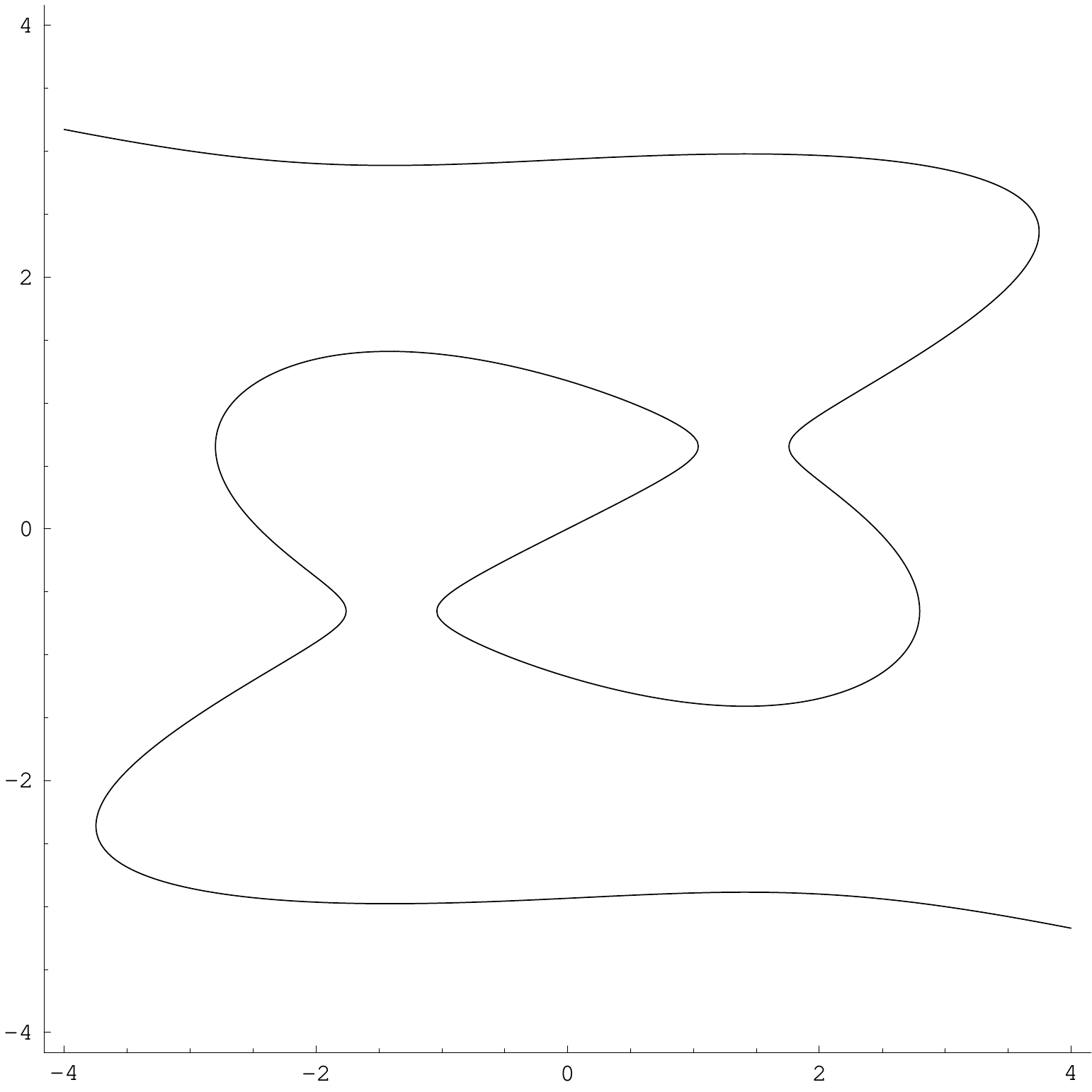}
\caption{Algebraic curve $P_{1/3}(x) = -P_{1/5}(y)$ defines the 1/3 to 1/5 map. Note the three nearly linear segments across the centre diagonal, matching Figure~11.}
\label{fig:12}
\end{center}
\end{figure}

It is worth noting that we can restrict the algebraic curve to the segments shown in Figure~12 by plotting only those points $(x,y)$ such that $P_{1/3}(x) = -P_{1/5}(y)$ and $|P_{1/3}(x)| \leq 4.$

\section{Horizontal similarity: Algebraic curves}

We have seen in Sections~7 and 8 that the horizontal maps appear as nearly linear sections of algebraic curves, in the form
$$ | P_{\theta}(x) = \pm P_{\theta'}(x') | \leq 4,$$
for the two characteristic polynomials $P_\theta$ and $P_{\theta'}$, where $\theta, \theta'$ are rational. The full algebraic curves are interesting in themselves, so in this section we present a few plots of the curves and point out some obvious patterns. 

\begin{figure}[t]
\begin{center}
$\begin{array}{ll }
   \includegraphics[width=2.5in]{Curves_1.pdf}   &  
   \includegraphics[width=2.5in,height=1.7in]{Curves_1_3.pdf}  \\
  \includegraphics[width=2.5in]{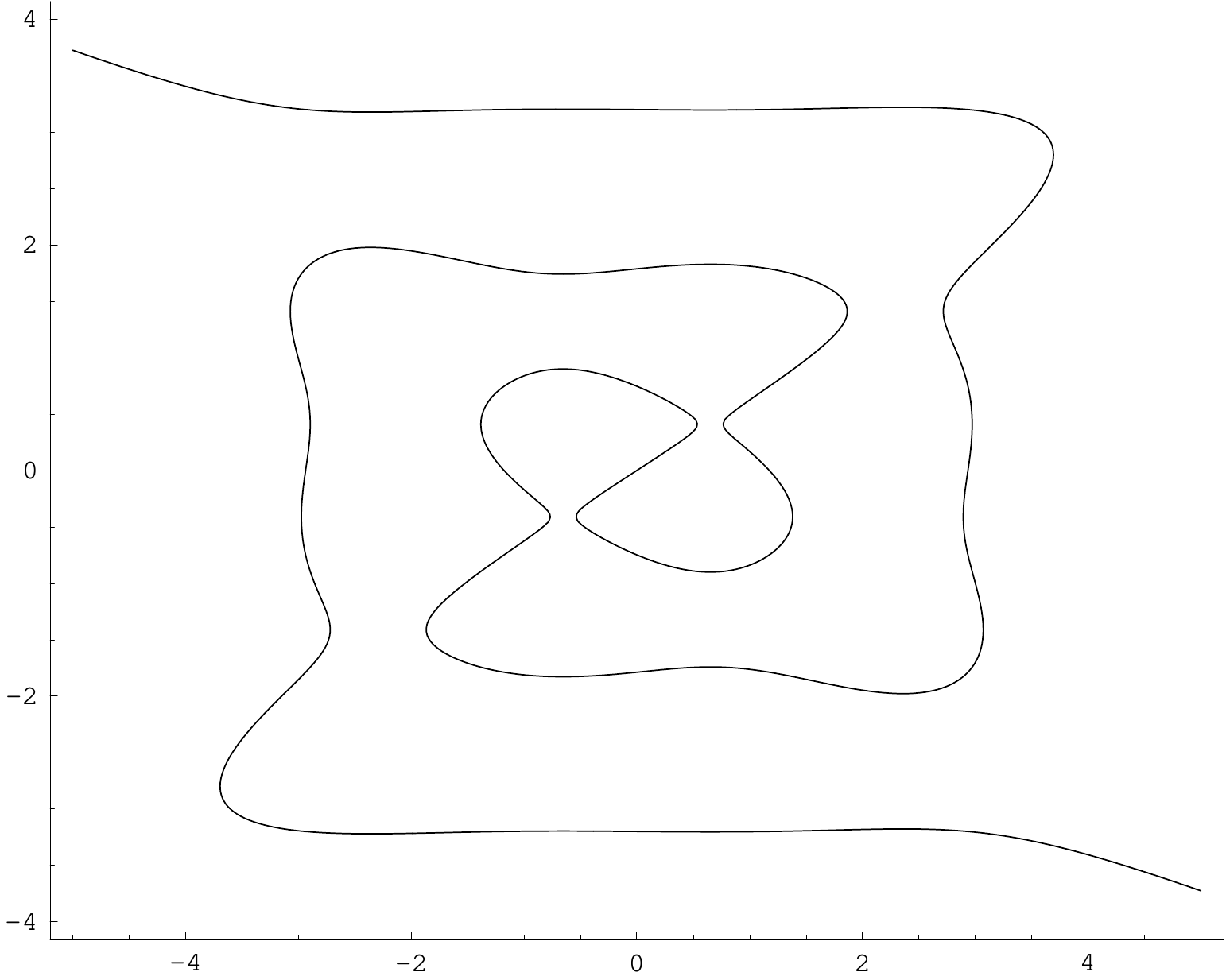}    &   
  \includegraphics[width=2.5in]{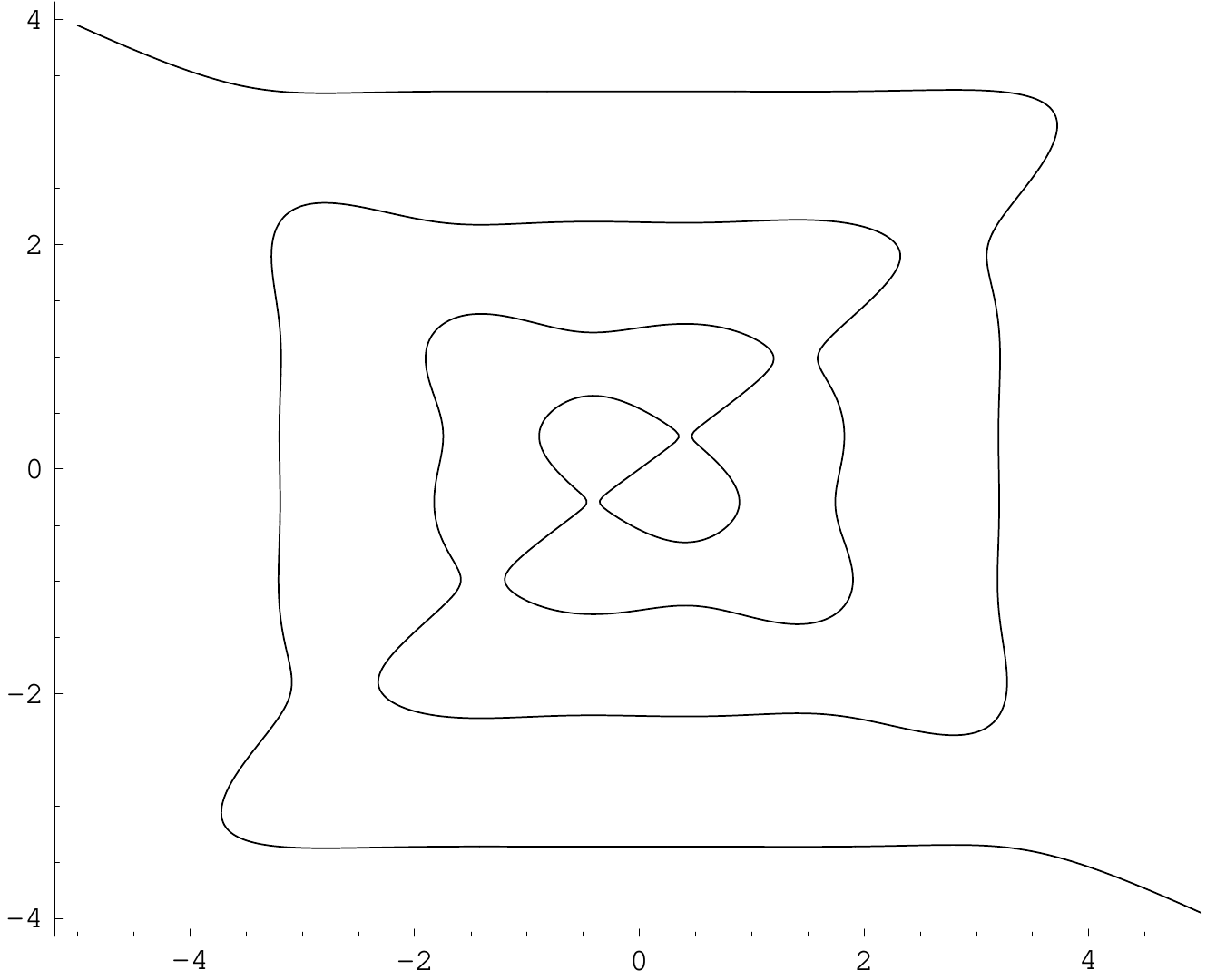} \\ 
  \includegraphics[width=2.5in]{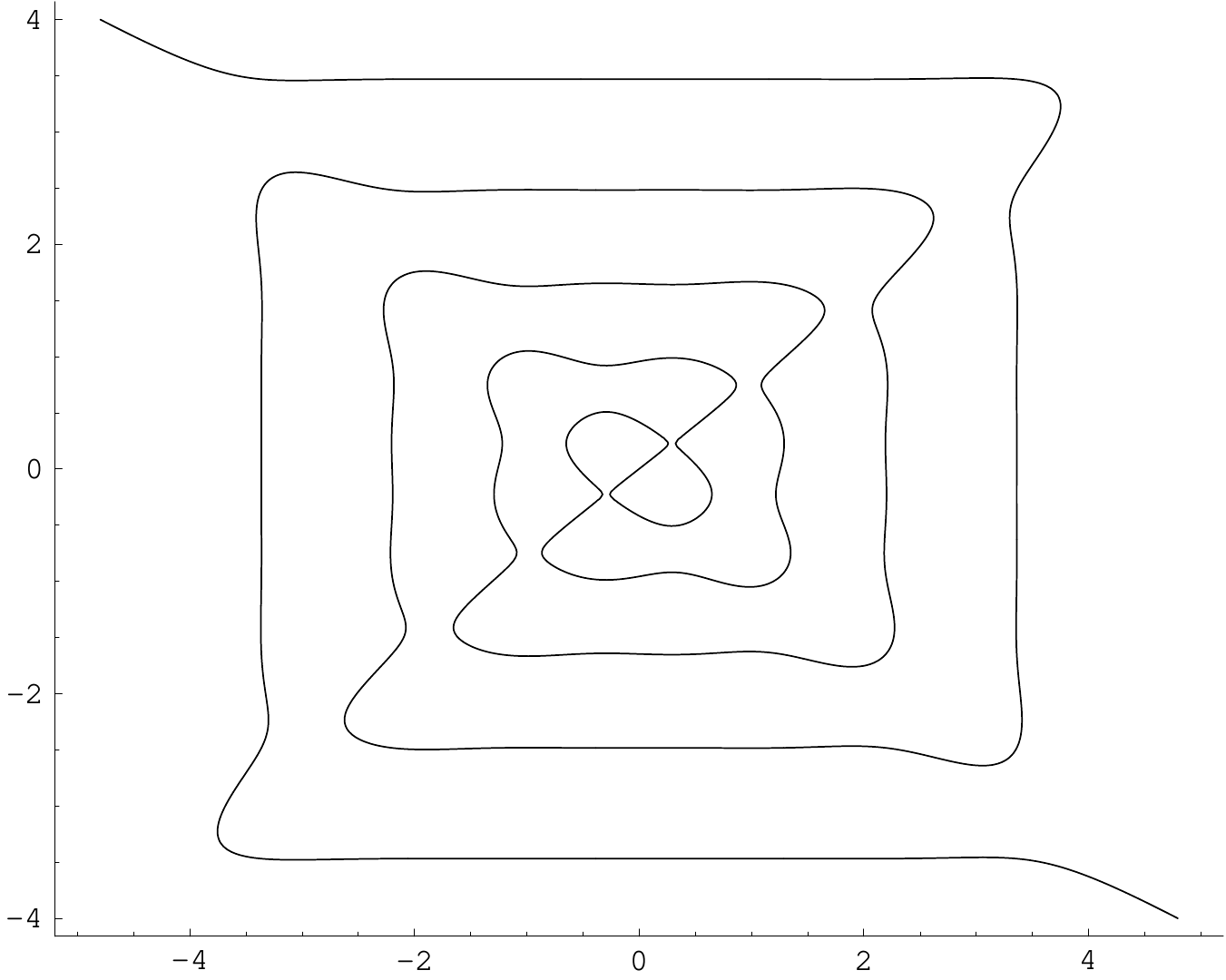} & 
\end{array}
$
\caption{A sequence of related algebraic curves, for rationals: $1\mapsto 1/3$,  $1/3 \mapsto 1/5$,  $1/5 \mapsto 1/7$,  $1/7 \mapsto 1/9$,  $1/9 \mapsto 1/11$. Notice the odd number of  nearly linear segments down the diagonal of each map.}
\label{mCurves_odds}
\end{center}
\end{figure}

Referring to Figure~\ref{mCurves_odds}, we see the algebraic curves $P_{\theta}(x) = P_{\theta'}(y)$ for rationals $\theta = 1/q \mapsto \theta' = 1/(q+2)$, with $q$ odd, have  a particularly simple graphs. There is a single connected component, that spirals around the origin, with more spirals as the denominator $q$ increases. We see $q$ nearly linear segments, which are the parts of the graph that define the relevant horizontal map on the fractal. We also see the odd symmetries of the curves, that point $(x,y)$ is in the curve if and only if $(-x,-y)$ is in it too.

\begin{figure}[t]
\begin{center}
$\begin{array}{ll }
   \includegraphics[width=2.5in]{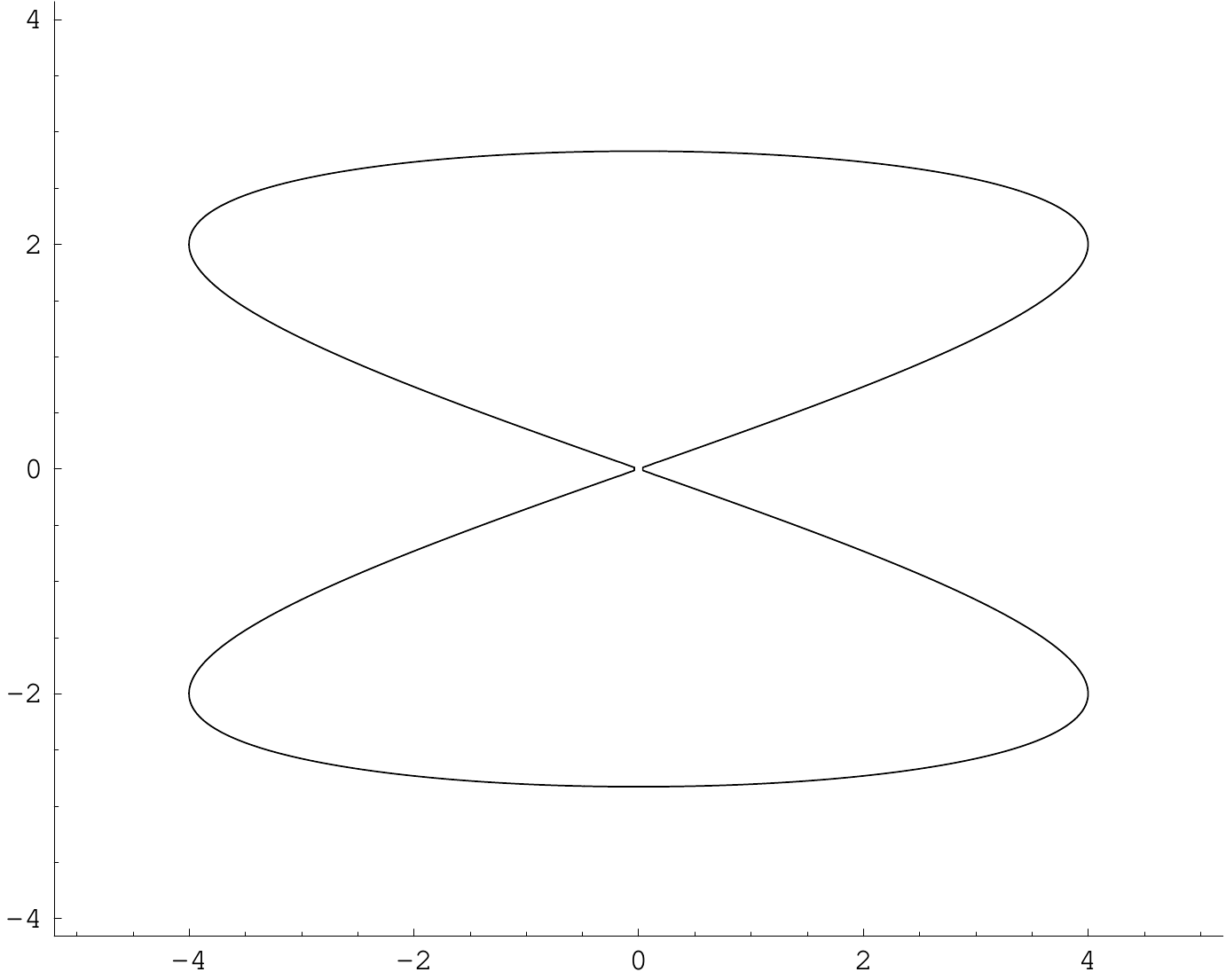}   &  
   \includegraphics[width=2.5in]{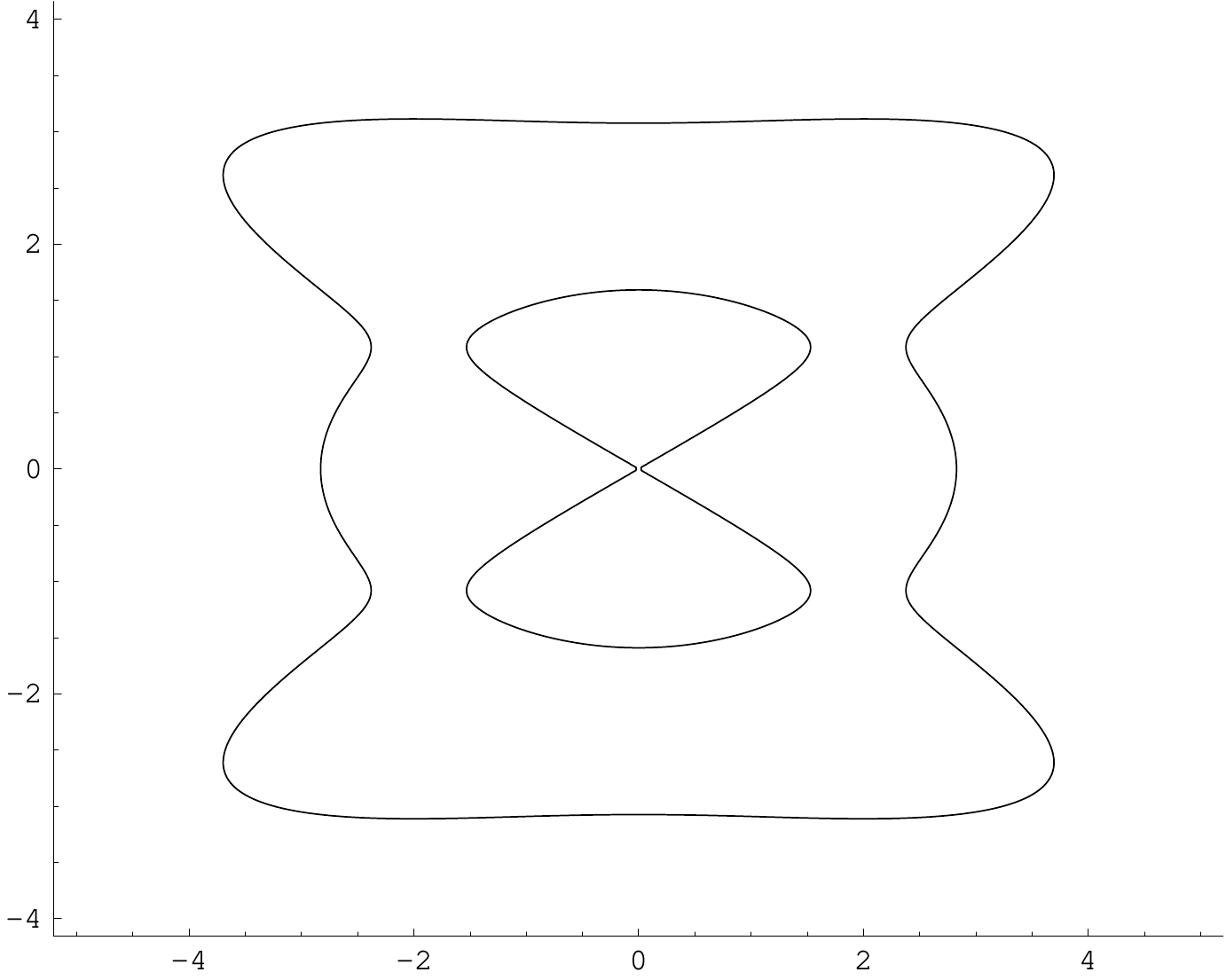}  \\
  \includegraphics[width=2.5in]{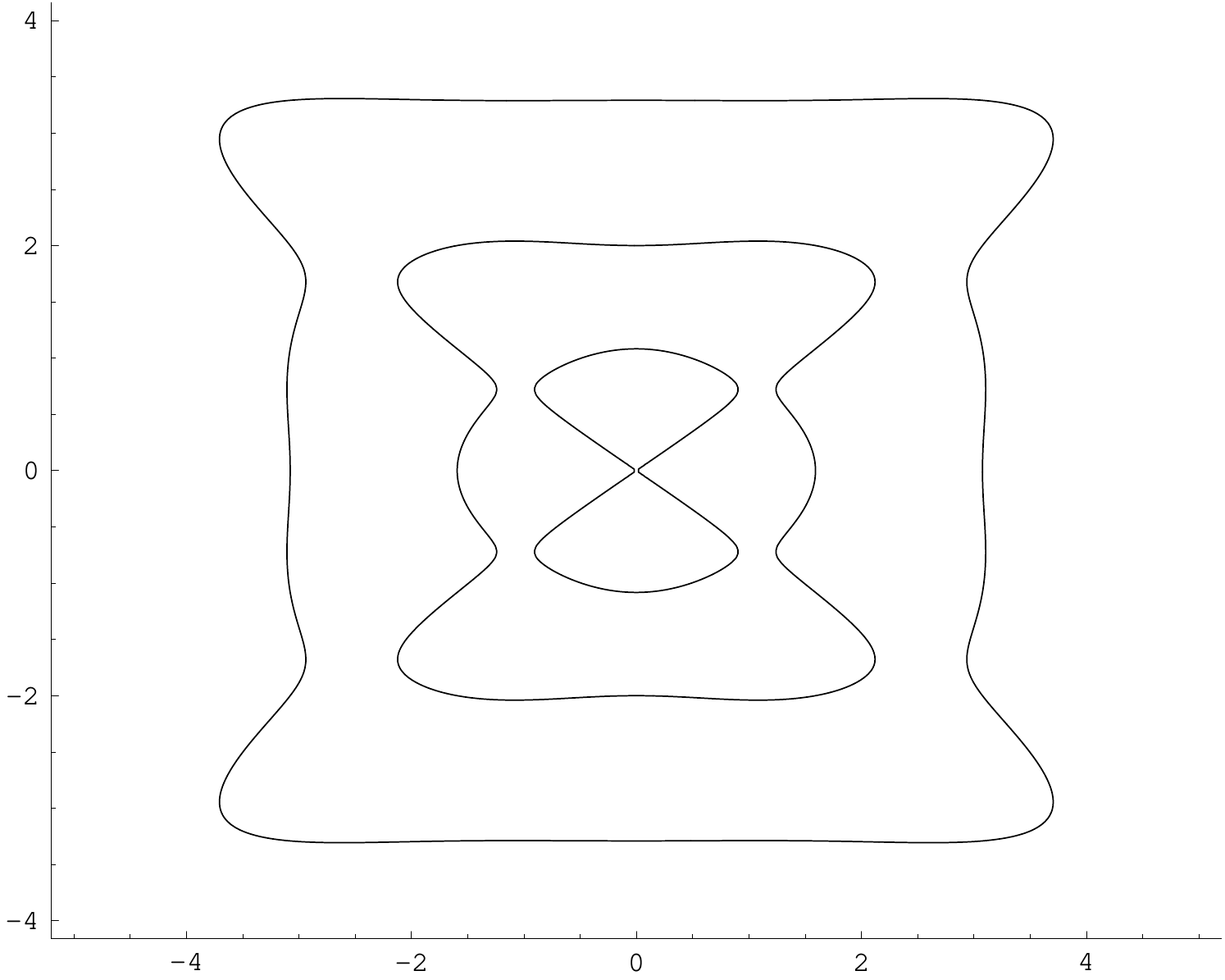}    &   
  \includegraphics[width=2.5in]{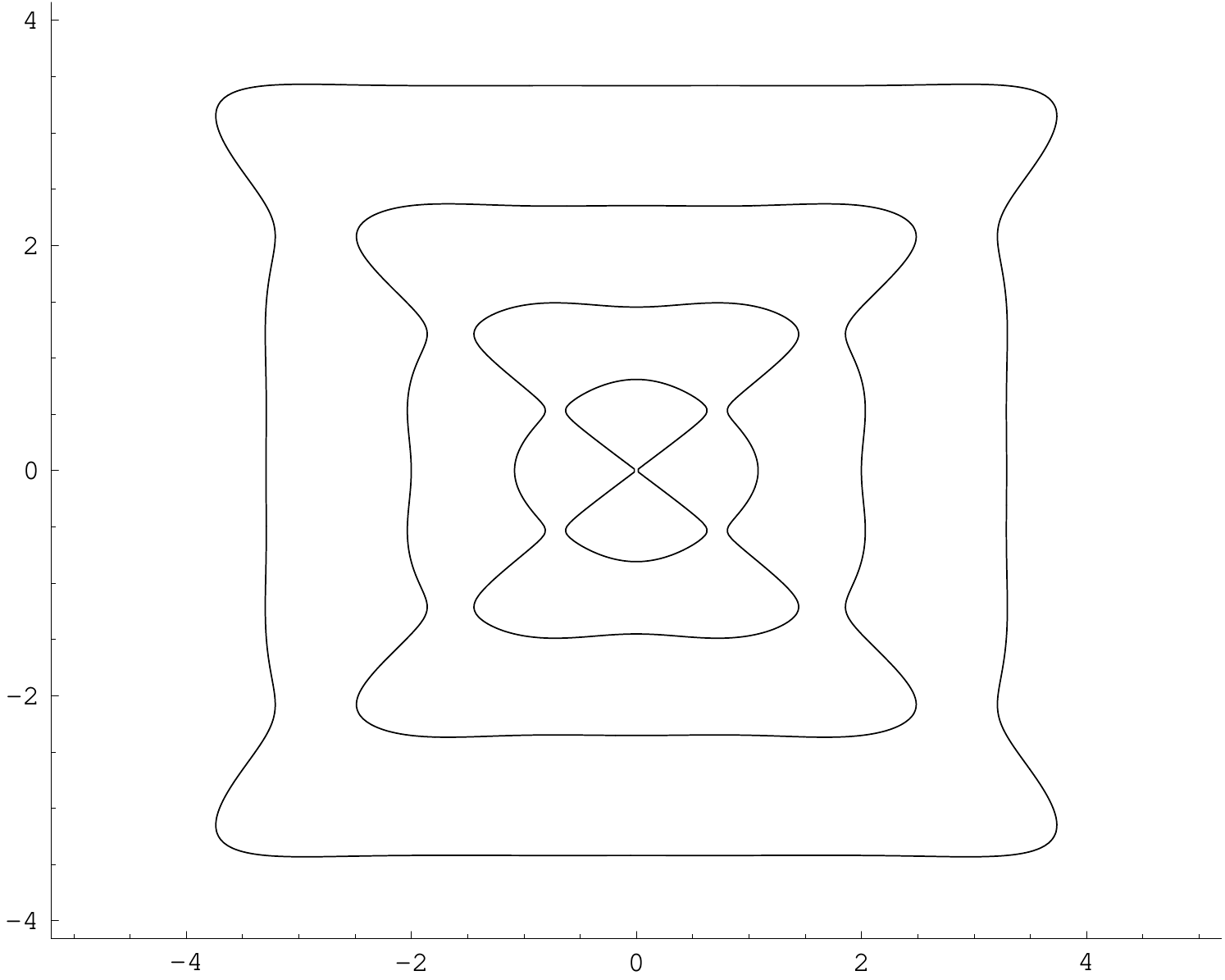} 
\end{array}
$
\caption{A sequence of related algebraic curves: $1/2\mapsto 1/4$,  $1/4 \mapsto 1/6$,  $1/6 \mapsto 1/8$,  $1/8 \mapsto 1/10$. Notice the symmetry, and the even number of  nearly linear segments down the diagonal of each map.}
\label{mCurves_evens}
\end{center}
\end{figure}

Referring to Figure~\ref{mCurves_evens}, we see the maps for rationals $1/q \mapsto 1/(q+2)$, with $q$ even, have  somewhat more complex graphs. There are $q/2$ connected components, that form something like figure eights around the origin, with more figure eights as the denominator $q$ increases. We see $q$ nearly linear segments, which are the parts of the graph that define the relevant horizontal map on the fractal. We also see the even symmetries of the curves, that point $(x,y)$ is in the curve if and only if the points $(\pm x, \pm y)$ are in it too.

\begin{figure}[ht]
\begin{center}
$\begin{array}{ll }
   \includegraphics[width=2.5in]{Curves_1_3.pdf}   &  
   \includegraphics[width=2.5in]{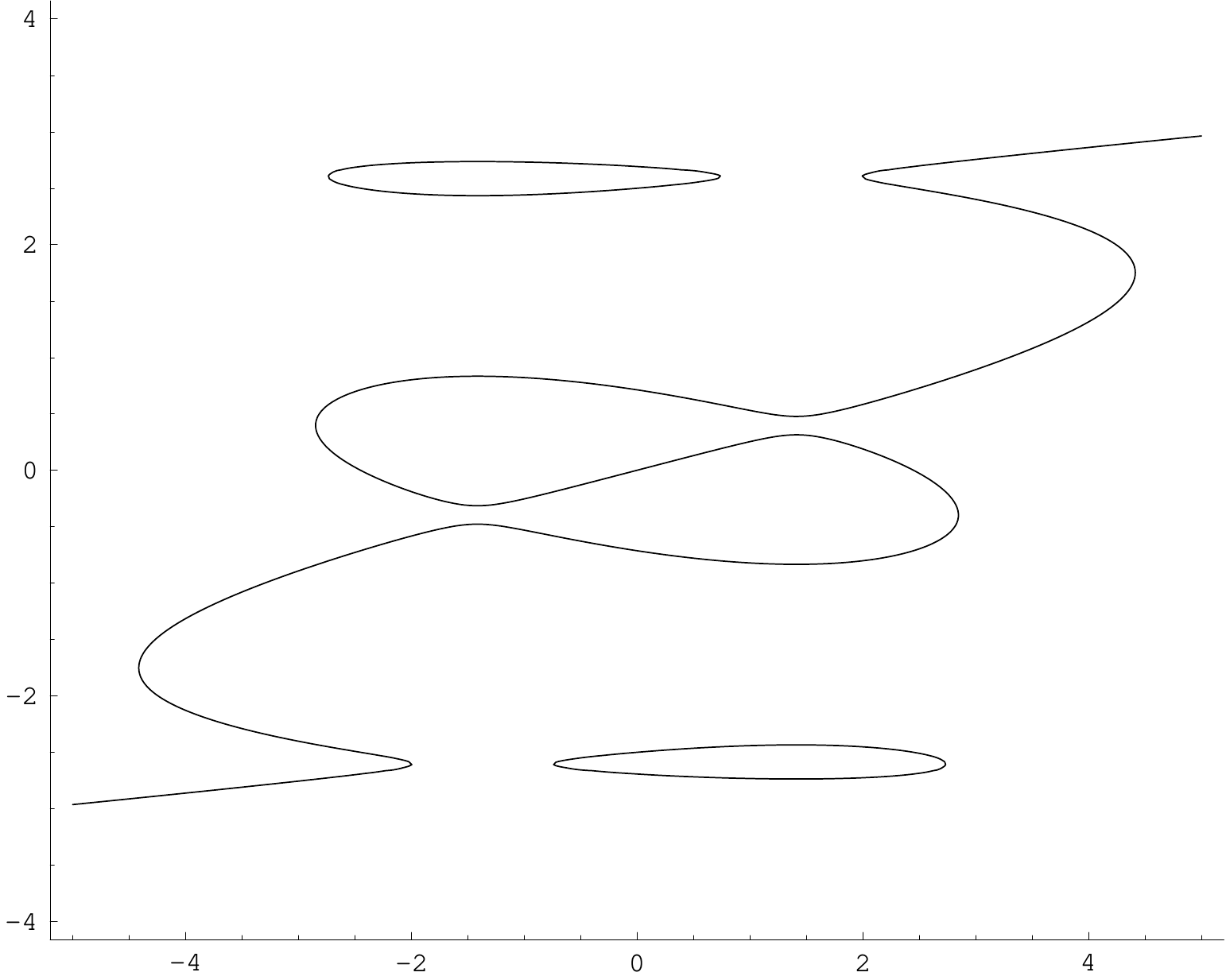}  
\end{array}
$
\caption{Algebraic curves with denominator 3: $1/3\mapsto 1/5$,  $2/3 \mapsto 2/7$. Notice the the three  nearly linear segments down the diagonal of each map, and the disconnected components on the right graph.}
\label{mCurves_threes}
\end{center}
\end{figure}

Figure~\ref{mCurves_threes} examines the curves for rational $\theta$ of the form $\theta=p/3$. In both cases ($p=1/3, 2/3$), there is odd symmetry, there are three nearly linear segments, and some spiraling of the curve about the origin. But for the case $\theta=2/3$, there are also two additional disconnected components to the curve, that seem to have nothing to do with the linear segments that are involved in the fractal map. 


In Appendix~3, we include more plots of these algebraic curves. We note with even denominator $q$, we have even symmetry in the plots, there are  $q$ nearly linear segments, and we see nested ``figure eights'' defining the key parts of the curve. For $p>1$,  there are disconnected components of the curve that seem to have nothing to do with the main part of the curve that defines the fractal map. 

We also observe the curves appear to get much more complex with increases in the numerator $p$, for the parameter $\theta = p/q$. There may be some connection between the various straight line segments in curves of different $\theta=p/q$ with the same denominator $q$, but these patterns are not entirely clear. 

We leave discussion of these curves and their connection with symmetries to future work.

\section{The similarity maps: General case}

The work in Sections~7 through 9  suggests an obvious candidate for the horizontal maps within spectral intervals: namely, a correspondence of points determined by two characteristic polynomials. Combining this observation with the earlier sections discussing the vertical maps and interval mapping, we can now specify the form of the general similarity maps of the butterfly, as follows:

Fix a matrix element $M$ in $ GL_2(\Z)$, represented as
\[ M = 
\left[
\begin{array}{rr}
a & b \\
c & d
\end{array}
\right].
\]
Fix an integer $r \geq 0$. The similarity map $S = S_{M,r,+}$ on the rational butterfly is described  as
\[ (x,\theta) \mapsto S(x,\theta) = (x',\theta'), \]
where the vertical component of the map is given by
\[ \theta \mapsto \theta' = \frac{ a\theta + b}{c\theta + d}, \]
while the horizontal component is given by the interval and polynomial correspondences. Specifically, 
 recall the spectrum at level $\theta = p/q$ is a set of $q$ intervals $I_k$, $k=1...q$, while at level $\theta'=p'/q'$ there are $q'$ intervals $I'_{k'}$, $k'=1...q'$. The map from $x$ at level $\theta$ to $x'$ at level $\theta'$ 
 is determined by the interval and polynomial conditions
\[ x\in I_k \mapsto x' \in I'_{r\cdot p' +k} =I_{k'} \mbox{ with } (-1)^{q+k}P_\theta(x) = (-1)^{q'+k'} P_{\theta'}(x'). \]
The sign choice for the polynomials is made so that $\pm P_\theta(x)$ and $ \pm P_{\theta'}(x')$ are both monotonically increasing on $I_k, I'_{k'}$, respectively. These two conditions determine the image point $x'$ uniquely, except possibly in one special case. Except for this case (described next), we have determined the similarity maps precisely. 

The special case to consider is when $x=0, \theta =p/q$ with $q$ even. In this case, the two intervals $I_{q/2}, I_{q/2 + 1}$ overlap at $x=0$ and may map to  two disjoint intervals $I'_{k'}$ and $I'_{k'+1}$. Here, we must decide where the point $x=0$ maps to, either an endpoint of $I'_{k'}$ or to an endpoint of $I'_{k'+1}$. It is convenient to chose to map to {\em both} points, making the function $S$ double-valued at this point $(x=0,\theta=p/q)$. We will see in the next section that the double-valued function is continuous.

There is another set of similarity maps, $S = S_{M,r,-}$, which differs from the above example by the way the interval maps are selected. For this maps, we chose
\[ x\in I_k \mapsto x' \in I'_{r\cdot (q'-p') +k} =I_{k'}, \]
which we observed in Section~6 as similarities that can occur.

Of course, not all such maps $S_{M,r,\pm }$ necessarily appear as similarities of the rational butterfly. For instance, from Theorem~1, we know the matrix $M$ must be chosen from the semigroup $G \subset GL_2(\Z)/{\pm I}$ described in Theorem~1. The integer $r$ must be small enough that
\[ r\cdot p' + q \leq q' \]
in the ``$+$'' case, and
\[ r\cdot (q'-p') + q \leq q' \]
 in the ``--" case, for all $\theta=p/q$. Equivalently, we have
 \[ r \leq \min_\theta \frac{c\theta + d-1}{a\theta +b} \]
in the ``$+$'' case, and
 \[ r \leq \min_\theta \frac{c\theta + d-1}{(c-a)\theta +d-b} \]
 in the ``--" case.
  But these are the only apparent restrictions on $M,r,\pm$.

\section{The similarity maps: Gap labelling}

The spectral gaps forming the butterfly are conveniently labeled using a Diophantine equation, where integer parameters $(s,t)$ are fixed, and  the $k$-th gap in the spectrum $\theta=p/q$ is identified using the formula
\[ k = t*p - s*q. \]
This indexing is known as gap labelling, as described in \cites{bel90, gold09, kam03, ypma07}, where the parameters $(s,t)$ are related to Chern numbers. These parameters give a convenient way of labelling the ``wings'' in the butterfly. Note there is a limited range on the indices: for $t>0$ we require 
\[ 0 \leq s \leq t-1, \]
while for $t<0$ we require
\[ t \leq s \leq -1. \]

It is easy to check that the similarity map $S=S_{M,r,\pm}$ maps labeled gaps to gaps, as stated in the following:

\begin{theorem}
The similarity map $S=S_{M,r,+}$ maps the gap labeled $(s,t)$ to the gap $(s',t')$ according to the formula
\[ 
\left[\begin{array}{c}s' \\t' \end{array}\right] =
(ad-bc)\left[\begin{array}{cc}a & b \\c & d\end{array}\right]
\left[\begin{array}{c}s \\t\end{array}\right] +
\left[\begin{array}{c}0\\r\end{array}\right] .
\]
Also, the similarity map $S=S_{M,r,-}$ maps the gap labeled $(s,t)$ to the gap $(s',t')$ according to the formula
\[ 
\left[\begin{array}{c}s' \\t'\end{array}\right] =
(ad-bc)\left[\begin{array}{cc}a & b \\c & d\end{array}\right]
\left[\begin{array}{c}s \\t\end{array}\right] +
\left[\begin{array}{c}r\\r\end{array}\right] .
\]
\end{theorem}

{\em Proof:} We consider $S=S_{M,r,+}$. Note the matrix $M$ in the similarity determines the map from angle $\theta = p/q$ to $\theta'=p'/q'$ according to the matrix formula
\[ 
\left[\begin{array}{c}p' \\q'\end{array}\right] =
\left[\begin{array}{cc}a & b \\c & d\end{array}\right]
\left[\begin{array}{c}p \\q\end{array}\right]  .
\]
In the case where the determinant $ad-bc$ is one, we invert to find
\[ p = (dp'-bq') \mbox{ and } q = (-cp' + aq'). \]

The index $r$ of the similarity tells us how intervals $I_k$ are mapped to intervals $I_{k'}$, where $k' = r\cdot p' + k$. Thus the $k$-th gap is mapped to gap $k' = r\cdot p' + k$. With gap $k=t*p-s*q$, we write
\begin{eqnarray}
k' &=& r\cdot p' + k \\
&=& r\cdot p' + t*p-s*q \\
&=& r\cdot p' + t(dp' - bq') - s(-cp'+aq') \\
&=& (cs + dt + r) p' - (as+bt) q' \\
&=& t'*p' - s'*q',
\end{eqnarray}
from which we read off $s' = (as+bt)$ and $t' = (cs +dt+ r)$, as desired.

The other cases (negative determinant, and $S=S_{M,r,-}$) are similar. \qed

\section{The similarity maps: Proof of continuity}

We show first continuity along the horizontal direction for the similarity maps described Section~10: 

\begin{theorem}
Given two rational parameters $\theta=p/q, \theta'=p'/q'$, and a correspondence between points in pairs of intervals $I_k, I'_{k'}$ in the spectra $Spec(h_\theta), Spec(h_{\theta'})$ respectively, given by
\[ x\in I_k \mapsto x' \in I_{k'}  = I'_{r\cdot p' +k}  \mbox{ with } (-1)^{q+k}P_\theta(x) = (-1)^{q'+k'} P_{\theta'}(x'), \]
then the correspondence is a continuous bijection on each interval. In particular, this yields a continuous map from the set of  points  $x \in Spec(h_\theta)$ into $x' \in Spec(h_{\theta'})$, with the possible exception at the point $x=0$, when $q$ is even.
\end{theorem}

{\em Proof:} The polynomial maps $(-1)^{q+k}P_\theta(x),  (-1)^{q'+k'} P_{\theta'}(x')$ are both continuous and monotonically increasing on the intervals in question, and thus have continuous inverses. So on each interval, we have a continuous bijection. For $q$ odd, the spectrum  $Spec(h_\theta)$ is a disjoint union of the intervals $I_1...I_q$, and so the collection extends to a continuous map into   $Spec(h_{\theta'})$. In the case where $q$ is even, two of the intervals $I_{q/2}$ and $I_{q/2+1}$ overlap at the point $x=0.$ So, except at this point, the correspondence extends to a continuous map on the union of intervals -- that is, it is continuous at $ Spec(h_\theta)\backslash 0$. At the point $x=0$, the correspondence maps the one point to (possibly) two points, one at the right endpoint of $I'_{r\cdot p' + q/2}$ and the other at the left endpoint of $I'_{r\cdot p' + q/2 +1}$. This double-valued splitting (although not a singled-valued function) is nevertheless continuous. \qed

\begin{figure}[t]
\begin{center}
\includegraphics[width=5in]{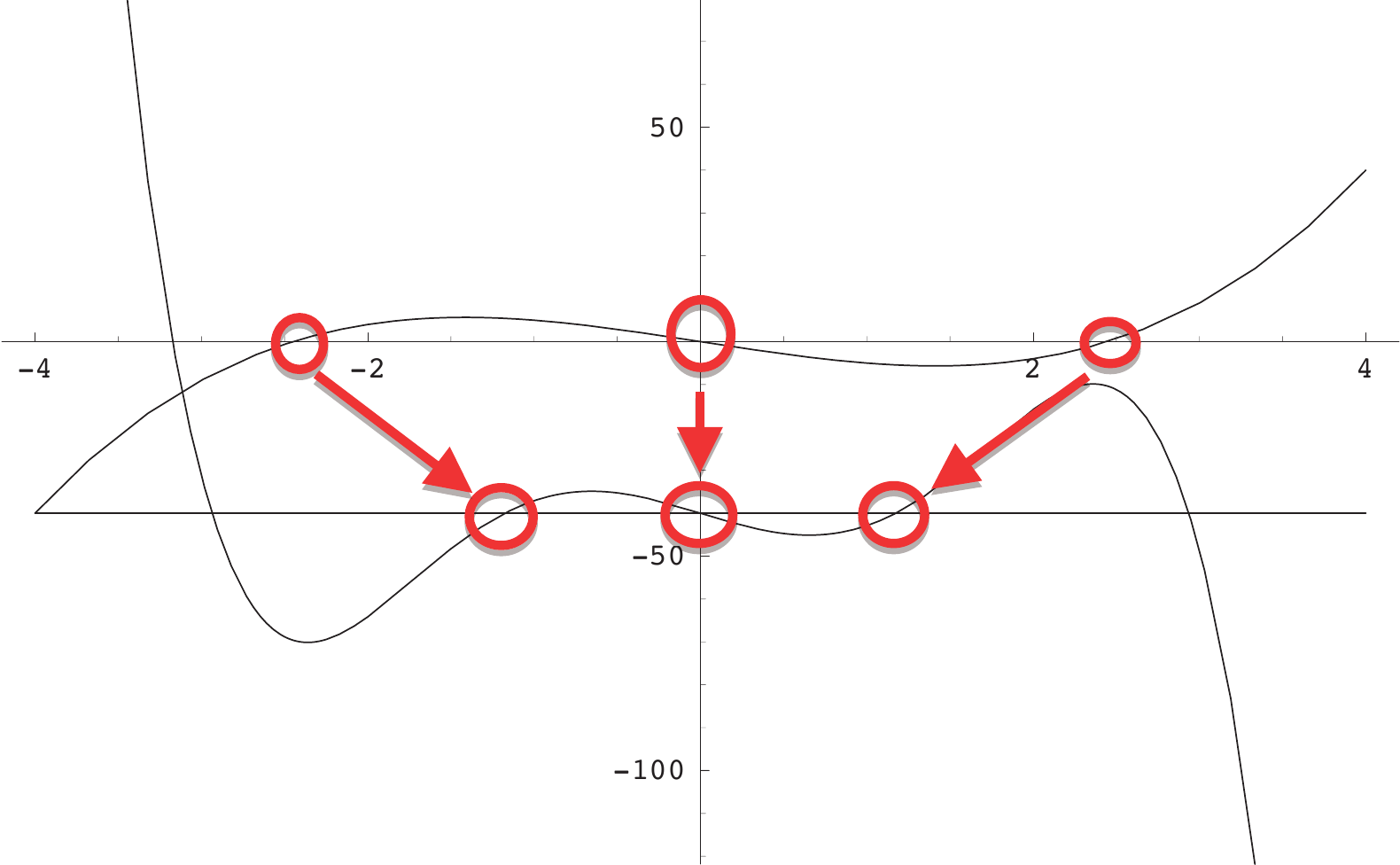}
\caption{Plot of cubic polynomial $P_{1/3}(x)$ and quintic $P_{1/5}(x')$. Match the zero crossing, move polynomials up and down, and observe a continuous correspondence of the zero crossings. This gives the continuous correspondence between line spectra.}
\label{fig:16}
\end{center}
\end{figure}

Figure~\ref{fig:16} indicates graphically how the bijection works on level $\theta=1/3$ to $\theta=1/5$.

We note that the correspondences of the form 
\[ x\in I_k \mapsto x' \in I_{k'} = I'_{r\cdot (q'-p') +k}  \mbox{ with } (-1)^{q+k}P_\theta(x) = (-1)^{q'+k'} P_{\theta'}(x'), \]
also are continuous, as in the theorem. This covers the similarity maps of the form $S_{M,r,-}$.

The numerical work suggests that piecing together these individual polynomial maps, at the different $\theta$ levels, gives a continuous map on the whole butterfly. We now state and prove this as a theorem:

\begin{theorem}
Suppose $S=S_{M,r,\pm}$ is a similarity of the rational butterfly, as described in Section~10. Then $S$ is a continuous (single-valued) map, except possibly at the points $(x=0,\theta=p/q)$ with $q$ even, where the map is double-valued and continuous.
\end{theorem}

{\em Proof:} We consider the case $S=S_{M,r,+}$. Fix a sequence of points $(x_n,\theta_n)$ converging to point $(x_*,\theta_*)$ and set \[  (x'_n,\theta'_n) = S(x_n,\theta_n). \]
By continuity of the linear fractional transformation,  the sequence
$\theta'_n$ converges to  $\theta'_* = (a\theta_* + b)/(c\theta_* + d)$. To show continuity of $S$, we only need to verify convergence of the $x'_n$. We will show every convergent subsequence of the image points  $(x'_n,\theta'_n)$ converges to the point $ (x'_*,\theta'_*) = S(x_*,\theta_*) $, which, by compactness, establishes continuity.

Take a convergent subsequence of the $(x'_n,\theta'_n)$ with limit $(x'_o,\theta'_*)$. Suppose in this subsequence, there are infinity many points where $\theta'_n = \theta'_*$ (and consequently, $\theta_n = \theta_*$). Applying Theorem~2, by continuity along the spectral lines at $\theta_*$, we can conclude that $x'_n$ converges to $x'_*$. (Except possibly at $x_*=0, \theta_*=p/q$, with $q$ even. Here, the map $S$ may be doubly valued, as discussed earlier.)

The other case to consider is when there are at most finitely  many $n$ with $\theta'_n = \theta'_*$. By restricting to the tail of the subsequence, and renumbering the subsequence, WLOG we may assume $\theta_n \neq \theta_*$ for all $n$. 

Set the notation consistently, so in the domain we have 
$\theta_n = p_n/q_n$, $\theta_* = p_*/q_*$, 
and points $x_n \in I_{k_n}$, $x_*\in I_{k_*}$, 
while in the range we have 
$\theta'_n = p'_n/q'_n$, $\theta'_* = p'_*/q'_*$, 
and points $x'_n \in I'_{k'_n}$, $x'_*\in I'_{k'_*}, x'_o\in I'_{k'_o}$. From the definition of $S$, we know $k'_n = r\cdot p'_n + k_n $ and $k'_* = r\cdot p'_* + k_*$. The value of $k'_o$ is to be determined (it will equal $k'_*$).

 We apply the trace $\tau_{\theta_n}$ to the spectral projection $\chi_{(-\infty, x_n)}(h_{\theta_n})$. By Lemma~7 in Appendix~2, we have
\[ \tau_{\theta_n}(\chi_{(-\infty, x_n)}(h_{\theta_n}))
= \frac{k_n-1}{q_n} + \frac{1}{q_n} F((-1)^{k_n+q_n}P_{\theta_n}(x_n)), \]
where the function $F$ is an integrated density of states on $h_o$.
By continuity of the trace and spectral projection (see~\cite{boca01}, Chap.~5), this quantity converges to
\[ \tau_{\theta_*}(\chi_{(-\infty, x_*)}(h_{\theta_*}))
= \frac{k_*-1}{q_*} + \frac{1}{q_*} F((-1)^{k_*+q_*}P_{\theta_*}(x_*)). \]

Since the sequence of ratios $p_n/q_n \neq p_*/q_*$ converges to $p_*/q_*$, we have that $q_n$ converges to infinity, so the terms above with $1/q_n$ in them converge to zero. Thus, from the limit of the traces, we conclude
\[ \lim_{n\rightarrow\infty} \frac{k_n}{q_n} = 
\frac{k_*-1}{q_*} + \frac{1}{q_*} F((-1)^{k_*+q_*}P_{\theta_*}(x_*)). \]

Similarly, by applying the trace to the sequence of spectral projections on the convergent sequence $x'_n \rightarrow x'_o$, we find 
\[ \tau_{\theta'_n}(\chi_{(-\infty, x'_n)}(h_{\theta'_n}))
= \frac{k'_n-1}{q'_n} + \frac{1}{q'_n} F((-1)^{k'_n+q'_n}P_{\theta'_n}(x'_n)), \]
converges to
\[ \tau_{\theta'_*}(\chi_{(-\infty, x'_o)}(h_{\theta'_*}))
= \frac{k'_o-1}{q'_*} + \frac{1}{q'_*} F((-1)^{k'_o+q'_*}P_{\theta'_*}(x'_o)). \]
Thus
\[ \lim_{n\rightarrow\infty} \frac{k'_n}{q'_n} = 
\frac{k'_o-1}{q'_*} + \frac{1}{q'_*} F((-1)^{k'_o +q'_*}P_{\theta'_*}(x'_o)). \]

With 
\[ \frac{k'_n}{q'_n} = \frac{r\cdot p'_n + k_n}{q'_n} 
= \frac{r\cdot(ap_n/q_n + b) + k_n/q_n}{cp_n/q_n + d}
\]
we take limits to obtain the identity
\[ \lim \frac{k'_n}{q'_n} = 
r\cdot \theta'_* + \frac{1}{c\theta_*+d}\lim \frac{k_n}{q_n}. \]
Replacing with the limiting values computed above, we obtain
\[ \frac{k'_o-1}{q'_*} + \frac{1}{q'_*} F((-1)^{k'_o+q'_*}P_{\theta'_*}(x'_o)) =
r\cdot \theta'_* + \frac{1}{c\theta_*+d}
\left( 
\frac{k_*-1}{q_*} + \frac{1}{q_*} F((-1)^{k_*+q_*}P_{\theta_*}(x_*))
\right).
\]
Multiplying through by $q'_*$ and simplifying yields
\[ \label{keyeq}
k'_o-1 + F((-1)^{k'_o+q'_*}P_{\theta'_*}(x'_o)) =
r\cdot p_*' + k_*-1 + F((-1)^{k_*+q_*}P_{\theta_*}(x_*)).
\]

Now, if $x_*$ is an interior point of the interval $I_{k_*}$, this last equation has $F$ taking a fractional value strictly between zero and one. Equating the integer part of the equation gives
\[ k'_o-1  =
r\cdot p_*' + k_*-1 ,
\]
from which we conclude $k'_o = k'_*$, which tells us that the limit point $x'_o$ and image point $x'_*$ are in the same spectral interval $I'_{k_*}$.

Equating the fractional part of the equation gives 
\[ F((-1)^{k'_*+q'_*}P_{\theta'_*}(x'_o)) =
 F((-1)^{k_*+q_*}P_{\theta_*}(x_*)),
\]
which, from the polynomial correspondence between $x_*$ and $x'_*$, yields 
\[ F((-1)^{k'_*+q'_*}P_{\theta'_*}(x'_o)) =
 F((-1)^{k'_*+q'_*}P_{\theta'_*}(x'_*)).
\]
Since $F$, composed with polynomial $(-1)^{k'_*+q'_*}P_{\theta'_*}(x)$ is strictly increasing on the interval $I'_{k_*}$, the equality implies $x'_o = x'_*$, as desired.

Now, when $x_*$ is at an endpoint of interval $I_k$, we use gap labeling to show the limit point $x'_o$ must be the corresponding endpoint of interval $I'_{k_*}$. Suppose $x_*$ is the left endpoint of interval $I_k$: in this case, equation~\ref{keyeq} reduces to 
\[ k'_o-1 + F((-1)^{k'_o+q'_*}P_{\theta'_*}(x'_o)) =
r\cdot p_*' + k_*-1 + 0.
\] 
so we know that $k'_o$ equals either $k'_*-1$ or $k'_*$. We just need to show it is equal to $k'_*$.

With the point $x_* \neq 0$ at a left endpoint, there must be a gap to the immediate left of $x_*$, which has some label $(s,t)$. 
By continuity of gap labeling~(Prop~11.11 in~\cite{boca01}),
the gap with this label forms an open set, so for $n$ sufficiently large, the point $x_n$ is to the right of this gap. Hence the interval $I_{k_n}$ is also to the right of the gap. Since the interval number must be bigger than the gap number, we have
\[ k_n > t*p_n - s*q_n. \]
The image point $x'_{k_n}$ is in interval $I'_{k_n}$ and its index satisfies
\begin{eqnarray}
k'_n &=& r\cdot p'_n + k_n \\
& > &  r\cdot p'_n + t*p_n - s*q_n \\
& = &  r\cdot p'_n + t*(dp'_n -b q'_n) - s*(-cp'_n + aq'_n) \\
&=&  (cs+dt+r)*p'_n - (as+bt)*q'_n \\
&=& t'*p'_n - s'*q'_n
\end{eqnarray}
where the indices $(s',t')$ label the image gap, as in Theorem~2. Thus the interval $I'_{k'_n}$ is to the right of the gap $(s',t')$ and hence so is the point $x'_n$. Since this gap is open, the limit point $x'_o$ is to the right of the gap, and thus the interval $I'_{k_o}$ containing point $x'_o$ is to the right. Consequently, we have
\[ k'_o >  t'p'_* - s'q'_* = k'_* - 1. \]
This eliminates the possibility that $k'_o = k'_* -1$, so we are left with $k'_o = k'_* $, from which we quickly conclude that the point $x'_o$ is the left endpoint of interval $I'_{k_*}$ and thus is equal to $x'_*$.

Handling a right endpoint is similar.  A similar analysis of the case $x_*=0$ at the common endpoint of two overlapping intervals shows the image sequence $x'_n$ could converge to endpoints of two different, disjoint intervals.
\qed

It would seem the result on the rational butterfly should extend by continuity to the full butterfly, including irrational values for vertical parameter $\theta$. We state this precisely as a conjecture.

\begin{conjecture}
Given a similarity $S=S_{M,r,\pm}$ of the rational butterfly, as described in Section~10, there is a unique continuous extension of $S$ to the full butterfly, possibly double-valued at certain points along the vertical line $x=0$. 
\end{conjecture}

We do not have a proof of this result. It seems a density argument and use of the continuous trace on the field of rotation algebras $A_\theta$ should lead to the result. 

It is curious to consider how the polynomial correspondences on the rational spectral lines might extend to irrational values of $\theta$, in which case there are no finite polynomials to describe the mapping.

\section{Three generators for the butterfly similarities}

There are three similarity maps of the butterfly which apparently generate all the similarities discussed above. There is the horizontal flip $H$ defined as the map
\[ H(x,\theta) = (-x,\theta). \]
There is the vertical flip $V$ defined as the map
\[ V(x,\theta) = (x, 1-\theta), \]
which can be represented in the form discussed in Section~10 as
\[ V = S_{B,0,+}, \]
with matrix $B = \left[ \begin{array}{cc}  -1 & 1 \\ 0 & 1  \end{array} \right]$. 
Finally, there is the similarity map discussed in Section~4, taking the butterfly to the bottom half given by 
\[ S(x, \theta) = (x', \frac{\theta}{\theta +1}), \]
where $x\in I_k$ in the k-th interval at level $\theta$ is mapped to $x'\in I'_k$ in the k-th interval at level $\theta'$. This map is represented as 
\[ S = S_{A,0,+},\]
for matrix $A = \left[ \begin{array}{cc}  1 & 0 \\ 1 & 1  \end{array} \right]$.

Some elementary calculations show how these three similarities combine algebraically. First, there are the two  obvious identities
\[ H^2 = I, \qquad V^2 = I, \]
and the commutation relation
\[ HV = VH. \]

The horizontal flip almost commutes with similarity $S$; in fact it just gives a shift in the interval indexing as follows:
\[HS_{A,0,+}H = S_{A,1,+}. \]
This shifting extends to powers of $S$, so we find
\[(HS_{A,0,+}H)^n = S_{A^n,n,+}. \]
By composing powers of $S$ and $HSH$, we obtain all similarities of the form
\[ S_{A^n,r,+} \qquad \mbox{ for $0\leq r \leq n$.} \]
We can also verify that such operators combine in the form
\[ (S_{A^{n'},r',+})( S_{A^n,r,+}) = S_{A^{n+n'},r+r',+}. \]

The vertical flip does not commute with similarity $S$, but instead introduces the similarity $S_{M,r,-}$ in the mix. Again, a straightforward calculation shows
\[ V S_{A^n,r,+} V = S_{BA^nB,r,-}. \] 

As noted in Theorem~1, the matrices $A,B$ generate a large semigroup of linear fractional transformations. From the calculations noted above, the three similarities of the butterfly generate enough similarities to cover this semigroup.

\begin{theorem}
The three similarities $H,V,S$ of the rational butterfly generate a semigroup of continuous (possibly double-valued) similarities of the butterfly, including maps of the form
\[ S_{M, r, +} \qquad \mbox{ and } \qquad S_{M,r,-} \]
where the matrices $M$ range over the semigroup of elements of $GL_2(\Z)/\pm I$ representing linear fractional transformations mapping the interval $[0,1]$ into itself. 
\end{theorem}

The double-valued character is, of course, from the interval splitting in the case of $\theta=p/q$ with $q$ even. There may be other similarities of the butterfly. The point is, we can at least see this large semigroup from $GL_2(\Z)$ appearing.

\section{Conclusions}

The Hofstadter butterfly, representing spectra of a continuous family of almost Mathieu operators, shows obvious fractal-like symmetry. By  investigating these symmetries numerically, we have catalogued the self-similarity maps of the butterfly using a semigroup of M\"{o}bius transformations in the vertical direction, indexed by elements in the matrix group $GL_2(\Z)$. This semigroup is generated by two matrices. In the horizontal direction, the self-similarity maps are given by algebraic curves determined by characteristic polynomials. Properties of the algebraic curves are demonstrated in a series of plots. These algebraic curves show nearly linear segments, demonstrating an almost linear behaviour in the horizontal component of the self-similarity maps. We proved continuity of the similarity maps on the rational butterfly, possibly double-valued at points with horizontal parameter $x=0$. The similarity maps are generated by exactly three continuous symmetries. We conjecture the semigroup of similarities on the rational butterfly extends to a family of continuous similarities on the full butterfly. 

\section*{Acknowledgments}

This work was supported in part by NSERC Discovery grants of the first two authors, and an NSERC Summer Research award of the third author. Numerical calculations were done in MATLAB~\cite{matlab} and rendered directly in PostScript, following a method similar to those discussed in~\cite{cass05}. Algebraic curves were rendered using Mathematica~\cite{mathca}.

We would like to thank the organizers  of the 2006 BIRS Workshop on  Operator Methods in Fractal Analysis, Wavelets, and Dynamical Systems, for encouraging us to present this work in an early form.


\nocite{avila06}
\nocite{bel82}
\nocite{bel90}
\nocite{brown64}
\nocite{cass05}
\nocite{choi90}
\nocite{connes94}
\nocite{gold09}
\nocite{hof76}  
\nocite{kam03} 
\nocite{lam07} 
\nocite{lam97}
\nocite{last94}
\nocite{matlab}
\nocite{puig03}
\nocite{ypma07}

\bibliography{FractalMathieu}

\pagebreak

\section*{Appendix 1}

{\bf Proof of Theorem 1}: We show the corresponding semigroup in $GL(2,\Z)/\{\pm I\}$ is generated by the two matrices $A$ and $B$. Suppose the matrix
\[ M = 
\left[
\begin{array}{rr}
a & b \\
c & d
\end{array}
\right]
\]
is in the semigroup. 
In the case that $b=0$, the determinant condition gives $a = \pm d = \pm 1$; the condition that the LFT maps $[0,1]$ into $[0,1]$ reduces WLOG to $a=d = 1$ and $c \geq 0$. This leave matrix $M$ in the form
\[ M = 
\left[
\begin{array}{rr}
1 & 0 \\
c & 1
\end{array}
\right] = A^c,
\]
with $c \geq 0$. Hence the matrix $M$ is generated by $A$ alone.

In the case that $b \neq 0$, WLOG (by mod-ing out by $\pm I$), we can assume $d \geq b >0$ as we know that $0$ maps to $b/d \in [0,1]$ under the LFT. By the determinant condition $ad - bc = \pm 1$,  we have that the gcd of $b,d$ is one, so we may apply the Euclidean algorithm to obtain a sequence of strictly positive quotients $q_0, q_1, q_2, \ldots$ and remainders $b=r_0, r_1, r_2, \ldots $ with 
\begin{eqnarray*}
d &=& q_0 r_0 + r_1 \\
b &=& q_1 r_1 + r_2 \\
r_1 &=& q_2 r_2 + r_3 \\
& \ldots & \\
r_{k-1} &=& q_k r_k + r_{k+1}, 
\end{eqnarray*}
where $r_k = gcd(b,d) = 1$ and $r_{k+1} = 0$. 
Applying these values to matrix $M$, we can factor as
\begin{eqnarray*}
\left[
\begin{array}{rr}
a & b \\
c & d
\end{array}
\right]
&=&
\left[
\begin{array}{rr}
0 & 1 \\
1 & q_0
\end{array}
\right]
\left[
\begin{array}{rr}
a_0 & r_1 \\
c_0 & r_0
\end{array}
\right] \\
&=&
\left[
\begin{array}{rr}
0 & 1 \\
1 & q_0
\end{array}
\right]
\left[
\begin{array}{rr}
0 & 1 \\
1 & q_1
\end{array}
\right]
\left[
\begin{array}{rr}
a_1 & r_2 \\
c_1 & r_1
\end{array}
\right] \\
&=&
\left[
\begin{array}{rr}
0 & 1 \\
1 & q_0
\end{array}
\right]
\left[
\begin{array}{rr}
0 & 1 \\
1 & q_1
\end{array}
\right]
\cdots
\left[
\begin{array}{rr}
0 & 1 \\
1 & q_k
\end{array}
\right]
\left[
\begin{array}{rr}
a_k & 0 \\
c_k & 1
\end{array}
\right],
\end{eqnarray*}
for some appropriate values of $a_0, a_1, \ldots, c_0, c_1, \ldots$.

Now, the matrix factors with the $q_j$ are generated by the $A,B$ matrices, since
\[ 
\left[
\begin{array}{rr}
0 & 1 \\
1 & q
\end{array}
\right]
= A^{q-1}BA.
\]
In the last factor,
\[
\left[
\begin{array}{rr}
a_k & 0 \\
c_k & 1
\end{array}
\right],
\] we know that $c_k$ is non-negative, else the LFT has a pole in $[0,1]$, which is not allowed. By the determinant condition, we have $a_k = \pm 1$. In the case $a_k = 1$, we have that the last factor appears as
\[
\left[
\begin{array}{rr}
1 & 0 \\
c_k & 1
\end{array}
\right] = A^{c_k}.
\]
In the case where $a_k = -1$, we may combine the last two factors to observe
\begin{eqnarray*}
\left[
\begin{array}{rr}
0 & 1 \\
1 & q
\end{array}
\right]
\left[
\begin{array}{rr}
-1 & 0 \\
c & 1
\end{array}
\right]
&=& 
\left[
\begin{array}{rr}
0 & 1 \\
1 & q
\end{array}
\right]
\left[
\begin{array}{rr}
1 & -1 \\
0 & 1
\end{array}
\right]
\left[
\begin{array}{rr}
c-1 & 1 \\
c & 1
\end{array}
\right] \\
&=&
\left[
\begin{array}{rr}
0 & 1 \\
1 & q-1
\end{array}
\right]
\left[
\begin{array}{rr}
c-1 & 1 \\
c & 1
\end{array}
\right] \\
&=& (A^{q-2}BA) (BA^c).
\end{eqnarray*}
For $q\geq 2$, we are done: these last factors appear as generated by matrices $A,B$, as desired.

In the case $q=1$, it easy to check that we must have more than one $q$-type matrix in the factorization, for if not, we get
\[ 
\left[
\begin{array}{rr}
a & b \\
c & d
\end{array}
\right]
=
\left[
\begin{array}{rr}
0 & 1 \\
1 & q_0
\end{array}
\right]
\left[
\begin{array}{rr}
-1 & 0 \\
c_0 & 1
\end{array}
\right], 
\]
where $q_0=1$, giving an LFT mapping 1 to $(c+1)/c$, which is outside the interval $[0,1]$. With a second  $q$-type matrix  in the factorization, we get
\[
M = \cdots (A^{q-1}BA)(A^{-1}BA)(BA^{c_k}) = \cdots (A^{q-1}B)(AA^{-1})(BA)(BA^{c_k})
 = \cdots A^q B A^{c_k}, \]
 which puts $M$ in the semigroup generated by $A,B$.
\qed

Here is a nice example, to show that the $a_k = -1$ case really does appear.
\[
\left[
\begin{array}{rr}
2 & 1 \\
3 & 2
\end{array}
\right]
=
\left[
\begin{array}{rr}
0 & 1 \\
1 & 2
\end{array}
\right]
\left[
\begin{array}{rr}
-1 & 0 \\
2 & 1
\end{array}
\right] 
\]
The matrix on the left gives an LFT which maps 0 to 1/2, and 1 to 3/5. This similarity map actually appears in the  Hofstadter 
butterfly, as we can see by examining Figure~1. We get $a_1 = -1$ in the factorization, which perhaps was not expected. Thus that little factor matrix on the right is not in the semigroup, since it maps 1 to -1/3, and thus does not map interval $[0,1]$ to itself.

However,  the original matrix on the left is actually in the semigroup, as shown in the proof.

\section*{Appendix 2}

We prove a basic result concerning the trace of certain spectral projections in the rotation algebra $A_\theta$. The machinery for this result is standard, and we borrow heavily from the results in \cite{boca01}.

\begin{lemma}
Let $\tau_\theta$ and $\tau_o$ be the (framed) tracial states on rotation algebras  $A_\theta, A_o$, respectively, with $\theta = p/q$.
Let $[a,b]$ be the k-th interval $I_k$ in the spectrum of the operator $h_\theta = u + u^* + v +v^*$ in algebra $A_\theta$. Then for $x\in [a,b]$,
\[ \tau_\theta(\chi_{[a,x)}(h_\theta)) = \frac{1}{q} \tau_o( \chi_{[-4,(-1)^{k+q}P_\theta(x))}(h_o)), \]
where $\chi_{[a,x)}$ is the spectral projection onto interval $[a,x)$. Consequently,
\[ \tau_\theta(\chi_{(-\infty,x)}(h_\theta)) = \frac{k-1}{q} + \frac{1}{q} \tau_o( \chi_{[-4,(-1)^{k+q}P_\theta(x))}(h_o)). \]
\end{lemma}

Note: the choice of the sign in front of polynomial $P_\theta(x)$ ensures that the map
\[ x \mapsto (-1)^{k+q}P_\theta(x) \]
is monotonic {\em increasing} on the interval $I_k$. In the earlier sections of this paper, it is convenient to write the trace-spectral projection function as the cumulative density $F(x)$ for the operator $h_o$, so we may write
\[ \tau_\theta(\chi_{(-\infty,x)}(h_\theta)) = \frac{k-1}{q} + \frac{1}{q} F((-1)^{k+q}P_\theta(x)), \]
where
\[ F(x) = \tau_o( \chi_{[-4,x)}(h_o)) \]
is a continuous, strictly increasing function mapping $[-4,4]$ onto $[0,1]$.

\vskip 5mm 
{\em Proof:} We do a direct calculation. The trace on $A_\theta$ is obtained from a representation of the algebra in $C(\T^2)\otimes M_q(\C)$ with generators
\[ u = \iota_1 \otimes U_o, v = \iota_2 \otimes V_o, \]
where the functions $\iota_1, \iota_2$ are the coordinate maps $\iota_1(z_1,z_2) = z_1,\iota_2(z_1,z_2) = z_2$, and the $q\times q$ matrices $U_o, V_o$ are the usual cyclic permutation and the diagonal with powers of $e^{2\pi i \theta}$. The trace in this frame is just
\[ \tau_\theta = \mu_2 \otimes tr_q, \]
where $\mu_2$ is the usual Haar measure on $\T^2$ and $tr_q$ is the normalized trace on $q\times q$ matrices. (See \cite{boca01}, page 11.)

Fix $x$ in $I_k$. The matrix-valued function $h = u+u^* + v + v^*$ has exactly one eigenvalue in the interval $I_k$, which moves continuously for parameter values $(z_1,z_2) \in \T^2$. This eigenvalue contributes to the trace precisely when it lies in the interval $[a,x]$. That is, we get a contribution precisely when the value of the polynomial $(-1)^{k+q}P_\theta(x)$ is smaller than $z_1^q + \overline{z_1}^q + z_2^q + \overline{z_2}^q$. We compute the trace in terms of a characteristic function on the torus $\T^2$, normalized by $q$, so 
\[
\tau_\theta( \chi_{[a,x)}(h_\theta))
= \frac{1}{q}
\int_{\T^2} \chi_{((-1)^{k+q}P_\theta(x) <  z_1^q + \overline{z_1}^q + z_2^q + \overline{z_2}^q)} \, d\mu_2.
\]
With $ y= (-1)^{k+q}P_\theta(x)$, and noting that the map $(z_1,z_2) \mapsto (z_1^q, z_2^q)$ has measure preserving inverse on the torus, we can write
\begin{eqnarray*}
\tau_\theta( \chi_{[a,x)}(h_\theta)) 
&=&
\frac{1}{q}
\int_{\T^2} \chi_{(y <  z_1^q + \overline{z_1}^q + z_2^q + \overline{z_2}^q)} \, d\mu_2 \\
&=&
\frac{1}{q}
\int_{\T^2} \chi_{(y <  z_1 + \overline{z_1} + z_2 + \overline{z_2})} \, d\mu_2 \\
&=&
\frac{1}{q}
\int_{\T^2} \chi_{(P_o(y) <  z_1 + \overline{z_1} + z_2 + \overline{z_2})} \, d\mu_2 \\
&=&
\frac{1}{q}\tau_o( \chi_{[-4,y)}(h_o) ).
\end{eqnarray*}
where we use the fact that $P_o(y) = y$ and so the last integral above yields the trace in $A_o$ applied to the spectral projection obtained by applying $\chi_{[-4,y]}$ onto $h_o$.

Using the relation $ y= (-1)^{k+q}P_\theta(x)$, we obtain
\[
\tau_\theta( \chi_{[a,x)}(h_\theta)) =
\frac{1}{q}\tau_o( \chi_{[-4,(-1)^{k+q}P_\theta(x))}(h_o) ),
\]
as desired. The trace evaluated on the full interval $(-\infty,x)$ picks up an additional contribution of $(k-1)/q$ from the previous $k-1$ intervals in the spectrum.
\qed

It is worth noting that this trace integral can be computed explicitly, the result involves an integral of $\arccos(y/2 + \cos(\theta))$.

\section*{Appendix 3}

We conclude with some final plots of the algebraic curves which determine the horizontal component of the similarity maps. Properties to notice are that the curves have even symmetry for denominator $q$ even, odd symmetry for $q$ odd, and in all cases there is a sequence of nearly linear segments along the diagonal line $y=x$. It is also notably that as the numerator $p$ increases, the algebraic curves become much more complicated.

\begin{figure}[ht]
\begin{center}
$\begin{array}{ll }
   \includegraphics[width=2.5in]{Curves_1_5.pdf}   &  
   \includegraphics[width=2.5in]{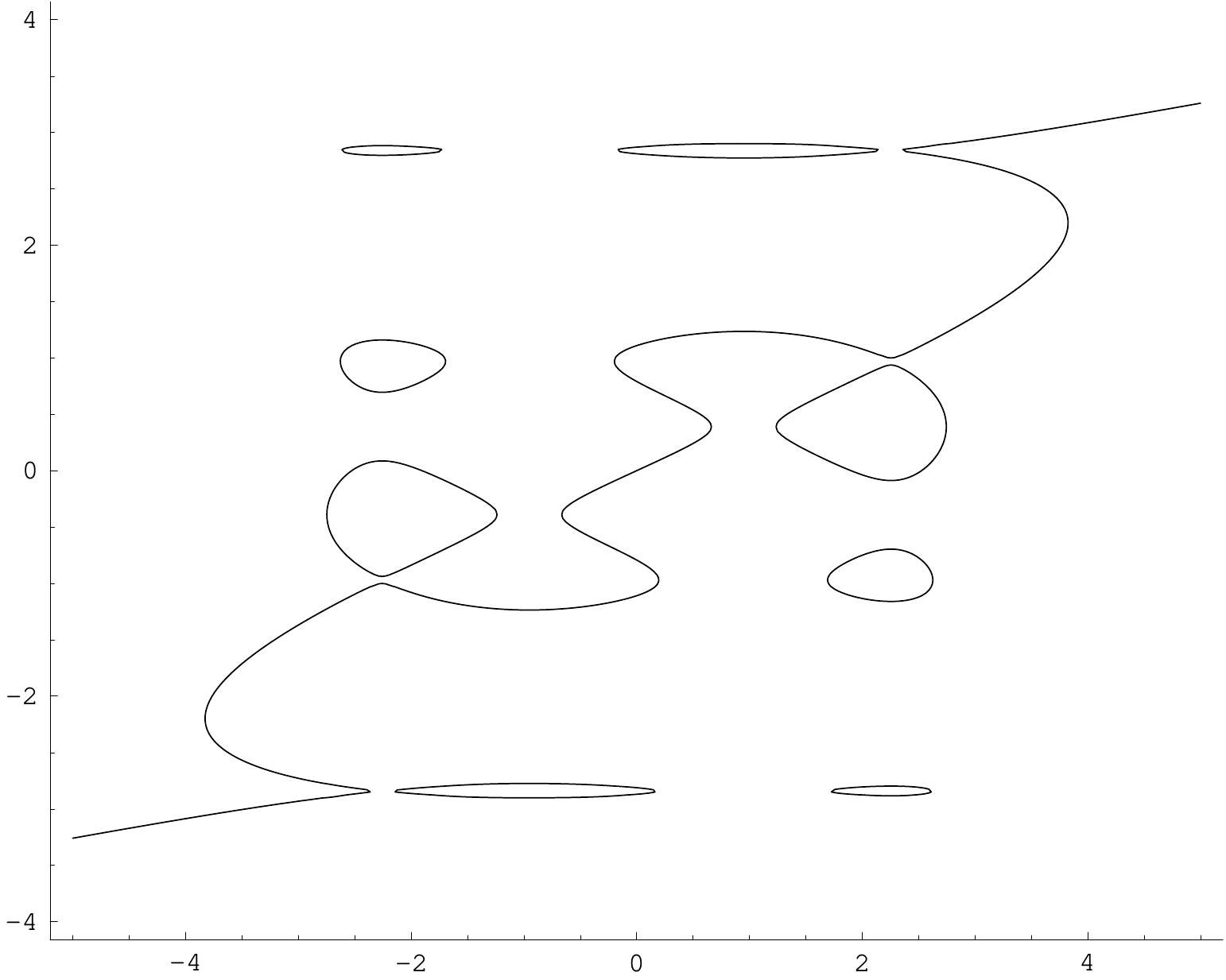}  \\
   \includegraphics[width=2.5in]{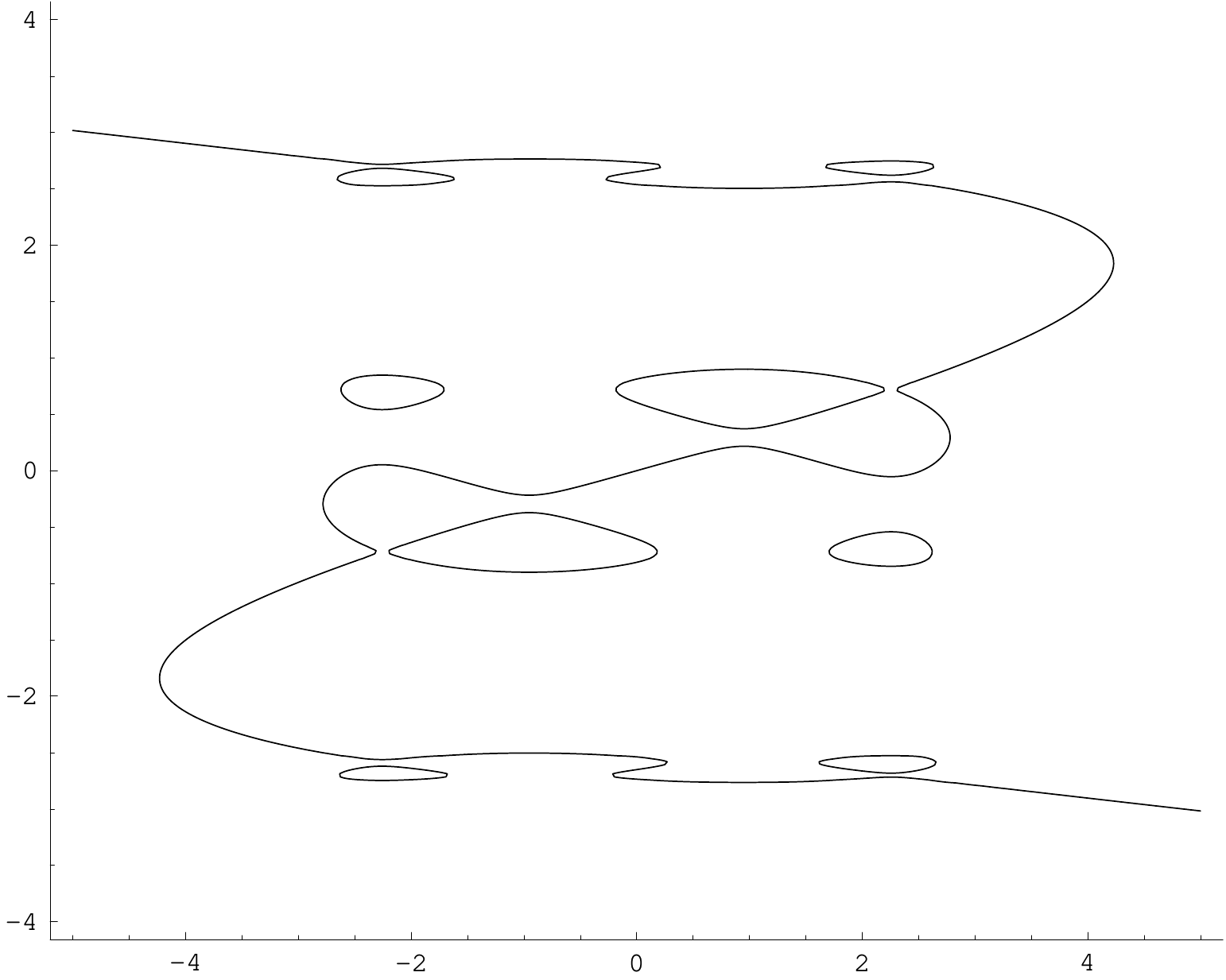}   &  
   \includegraphics[width=2.5in]{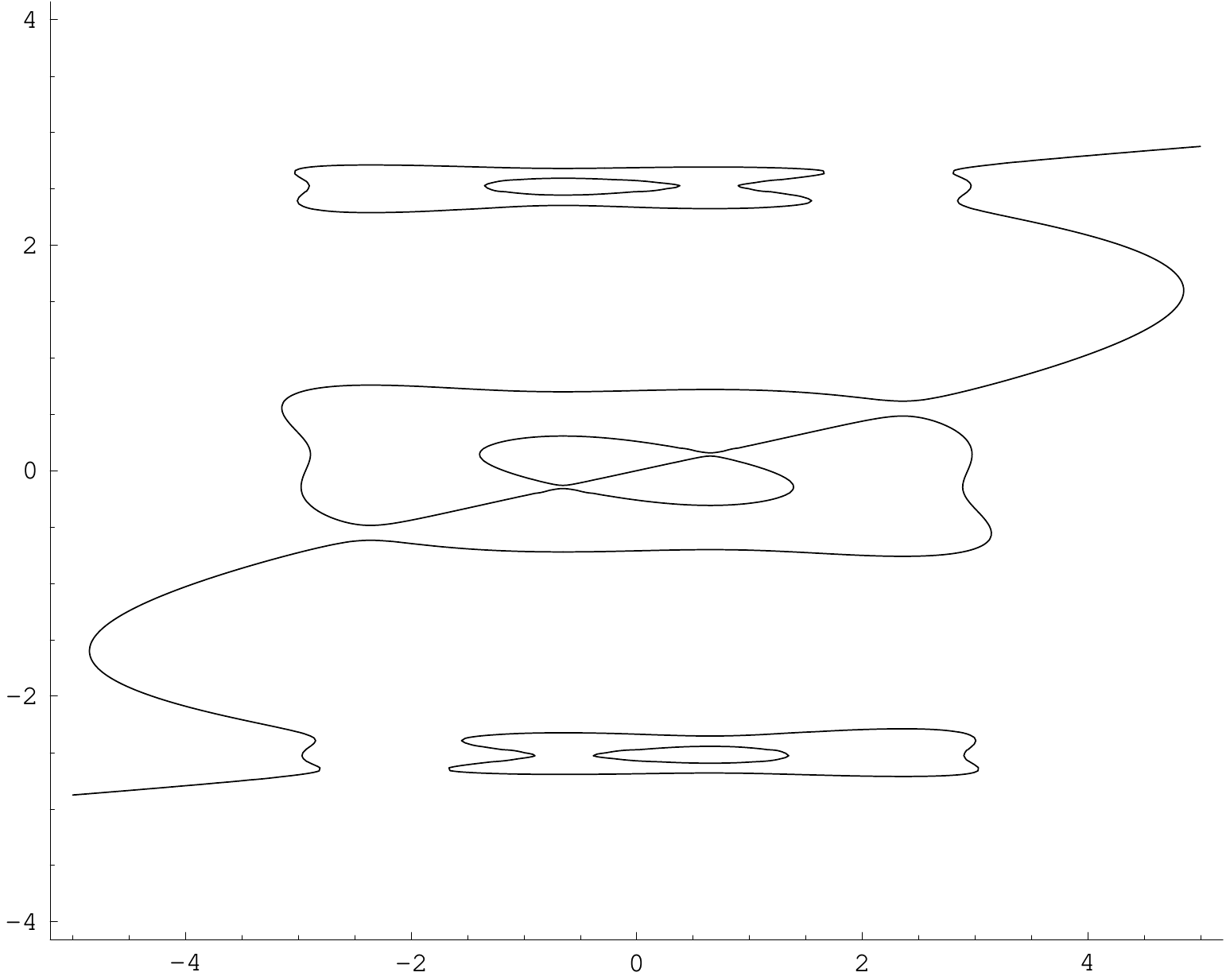}  \\
\end{array}
$
\caption{Algebraic curves with denominator 5: $1/5\mapsto 1/7$,  $2/5 \mapsto 2/9$, $3/5 \mapsto 3/11$, $4/5 \mapsto 4/13$. In each graph, there are five nearly linear segments down the diagonal of the graph. All but the first graph have some disconnected components.}
\label{mCurves_fives}
\end{center}
\end{figure}

\begin{figure}[ht]
\begin{center}
$\begin{array}{ll }
   \includegraphics[width=2.5in]{Curves_1_7.pdf}   &  
   \includegraphics[width=2.5in]{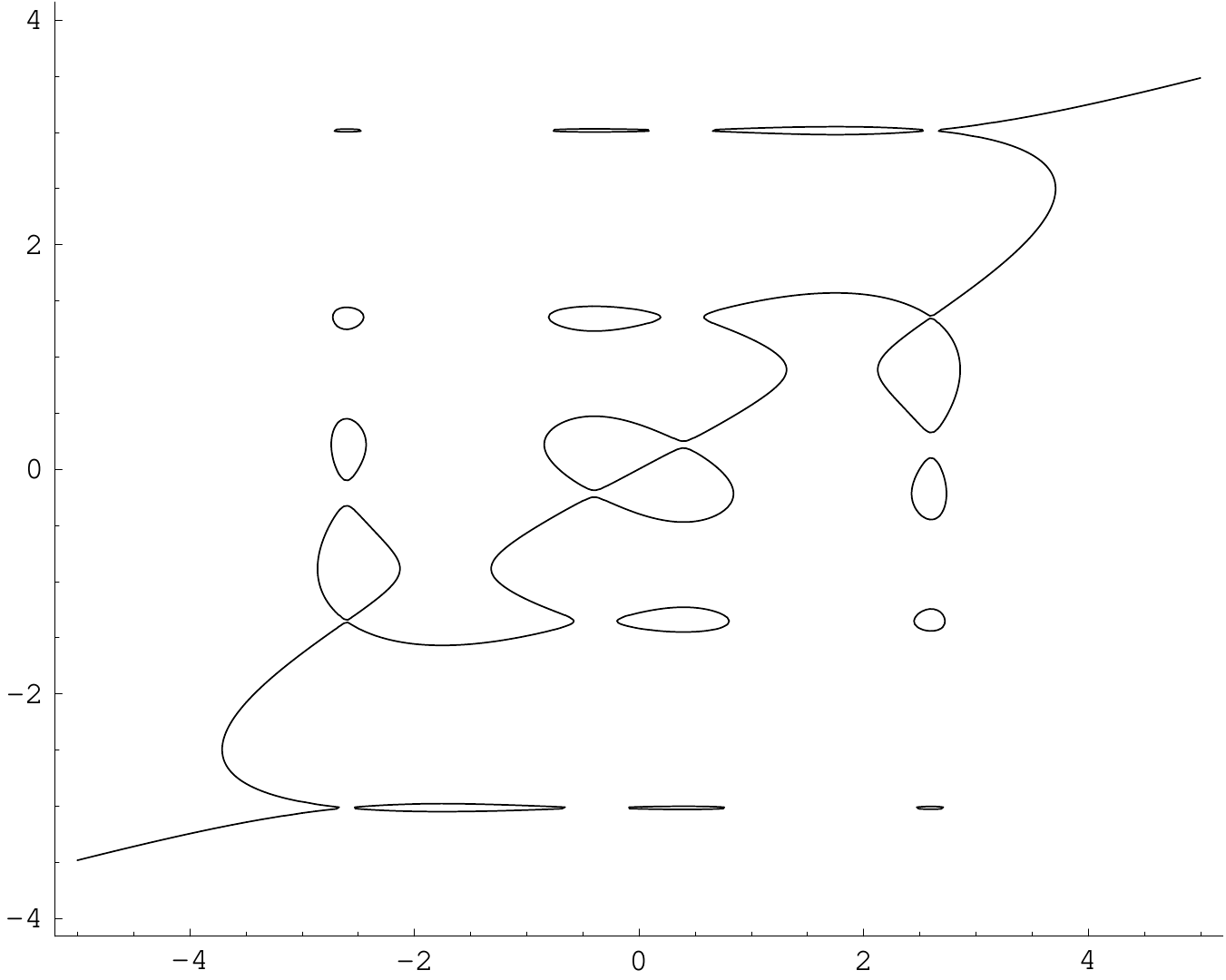}  \\
   \includegraphics[width=2.5in]{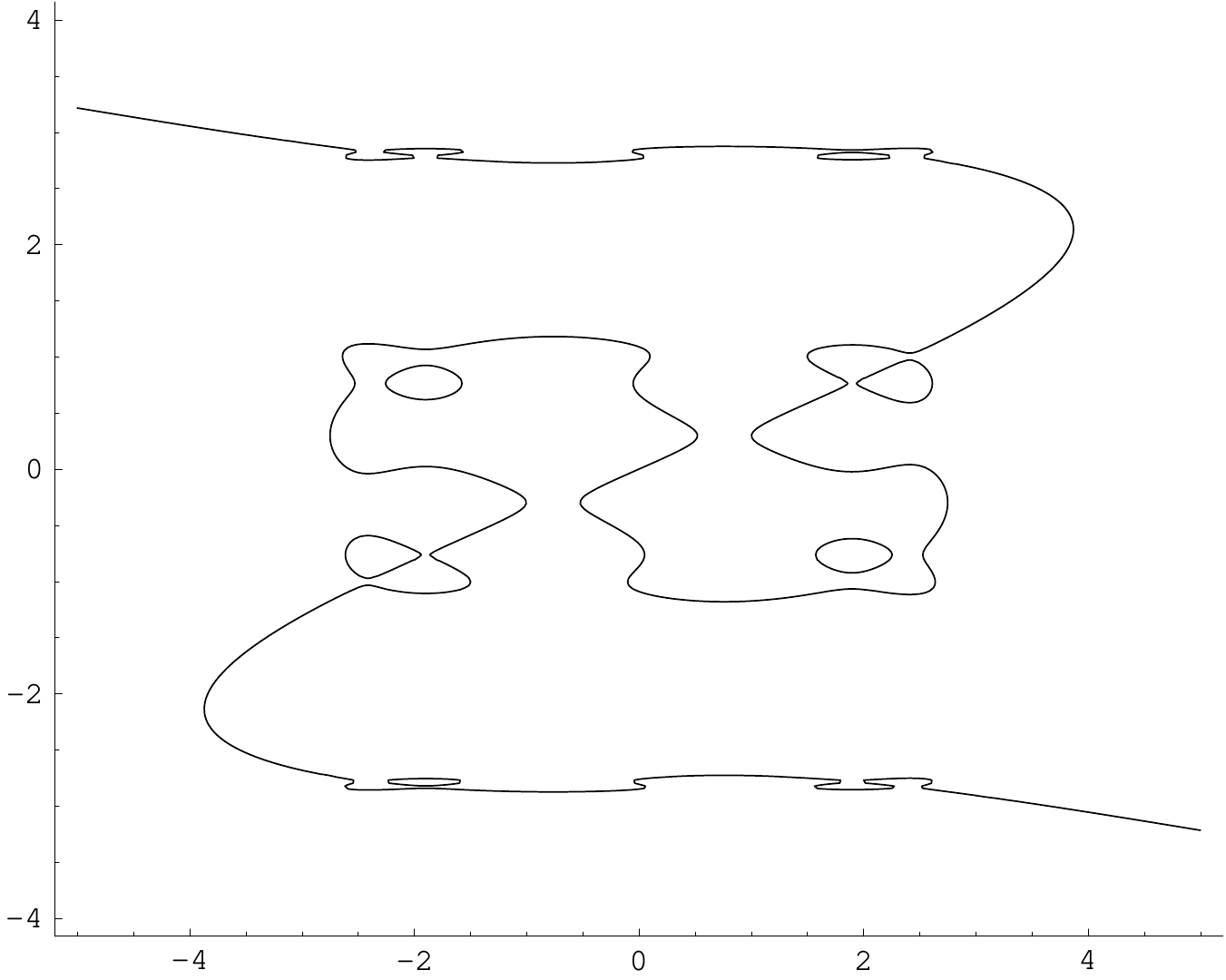}   &  
   \includegraphics[width=2.5in]{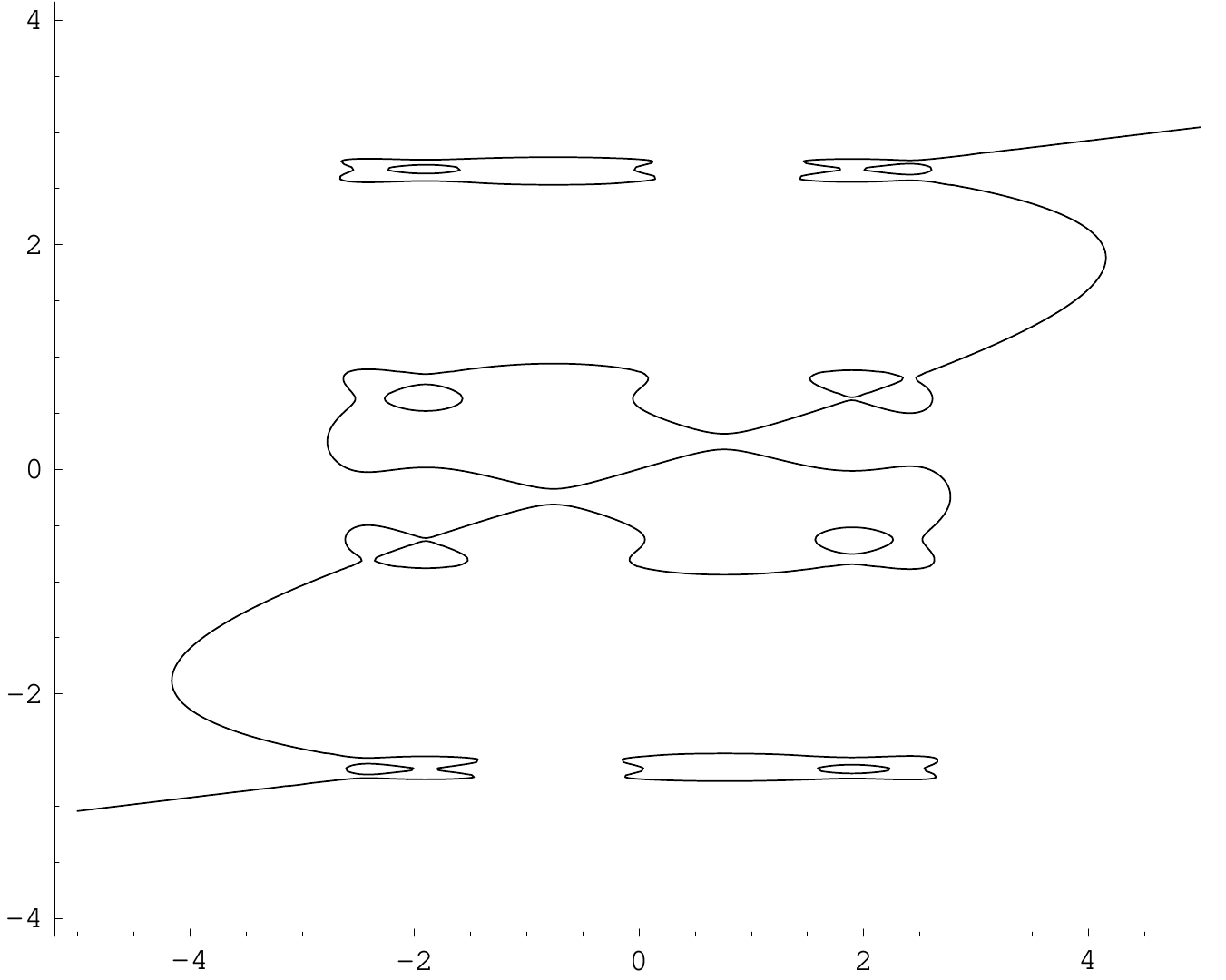}  \\
   \includegraphics[width=2.5in]{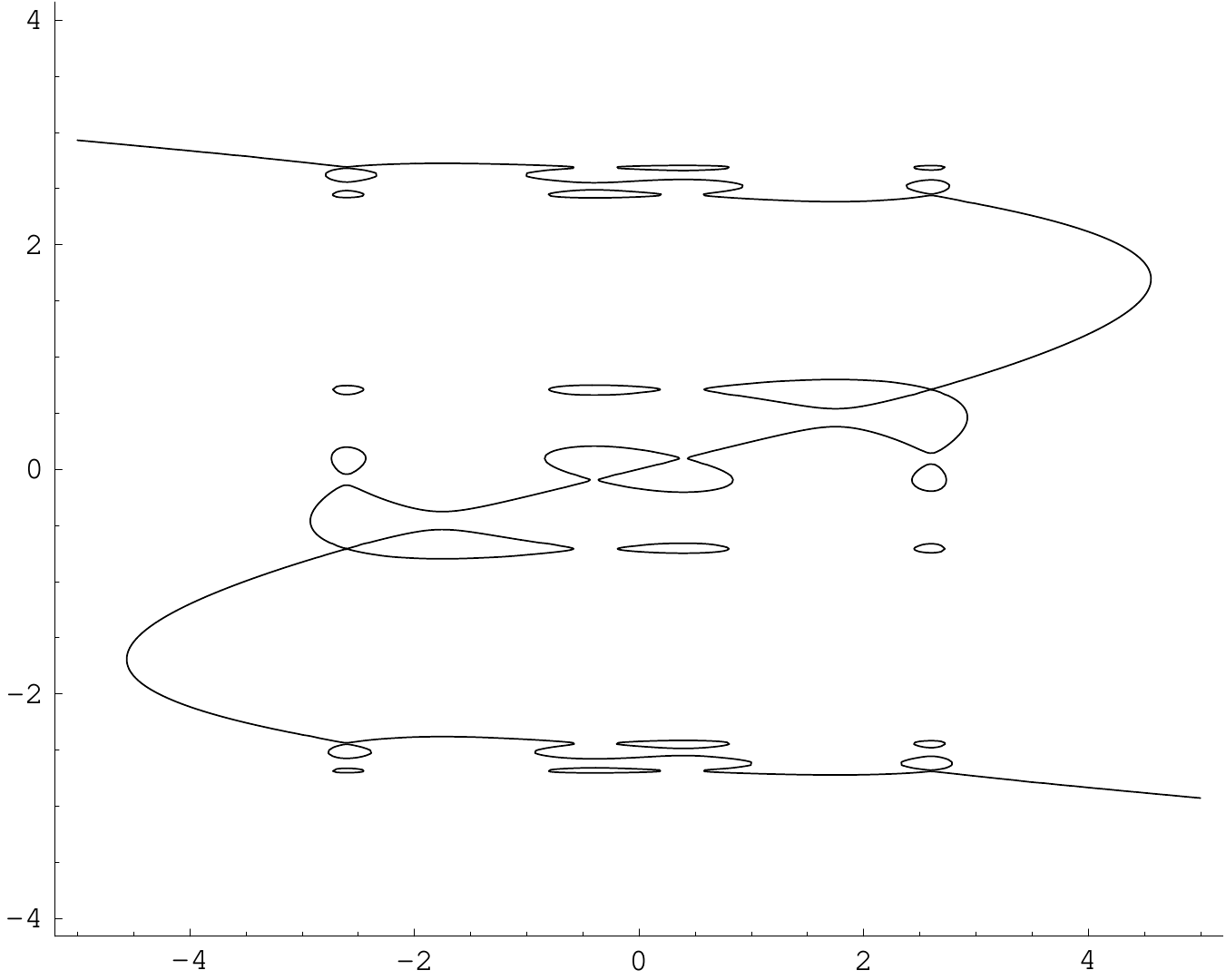}   &  
   \includegraphics[width=2.5in]{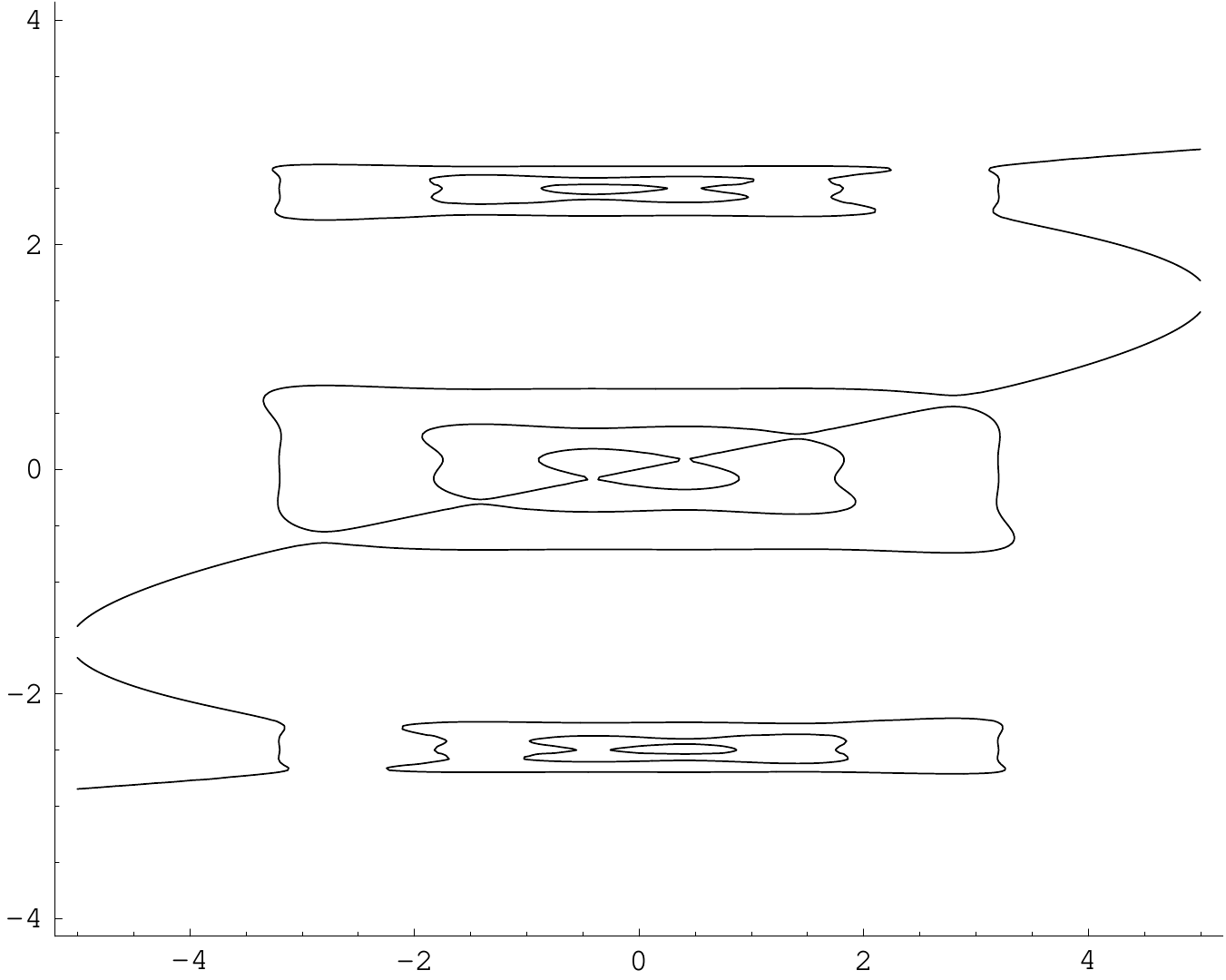}  \\
\end{array}
$
\caption{Algebraic curves with denominator 7: $1/7\mapsto 1/9$,  $2/7 \mapsto 2/11$, $3/7 \mapsto 3/13$, $4/7 \mapsto 4/15$, $5/7 \mapsto 5/17$, $6/7 \mapsto 4/19$. In each graph, there are seven nearly linear segments down the diagonal of the graph. All but the first graph have some disconnected components.}
\label{mCurves_sevens}
\end{center}
\end{figure}

.. \hfill ..
\pagebreak

\begin{figure}[h]
\begin{center}
\includegraphics[width=2.5in]{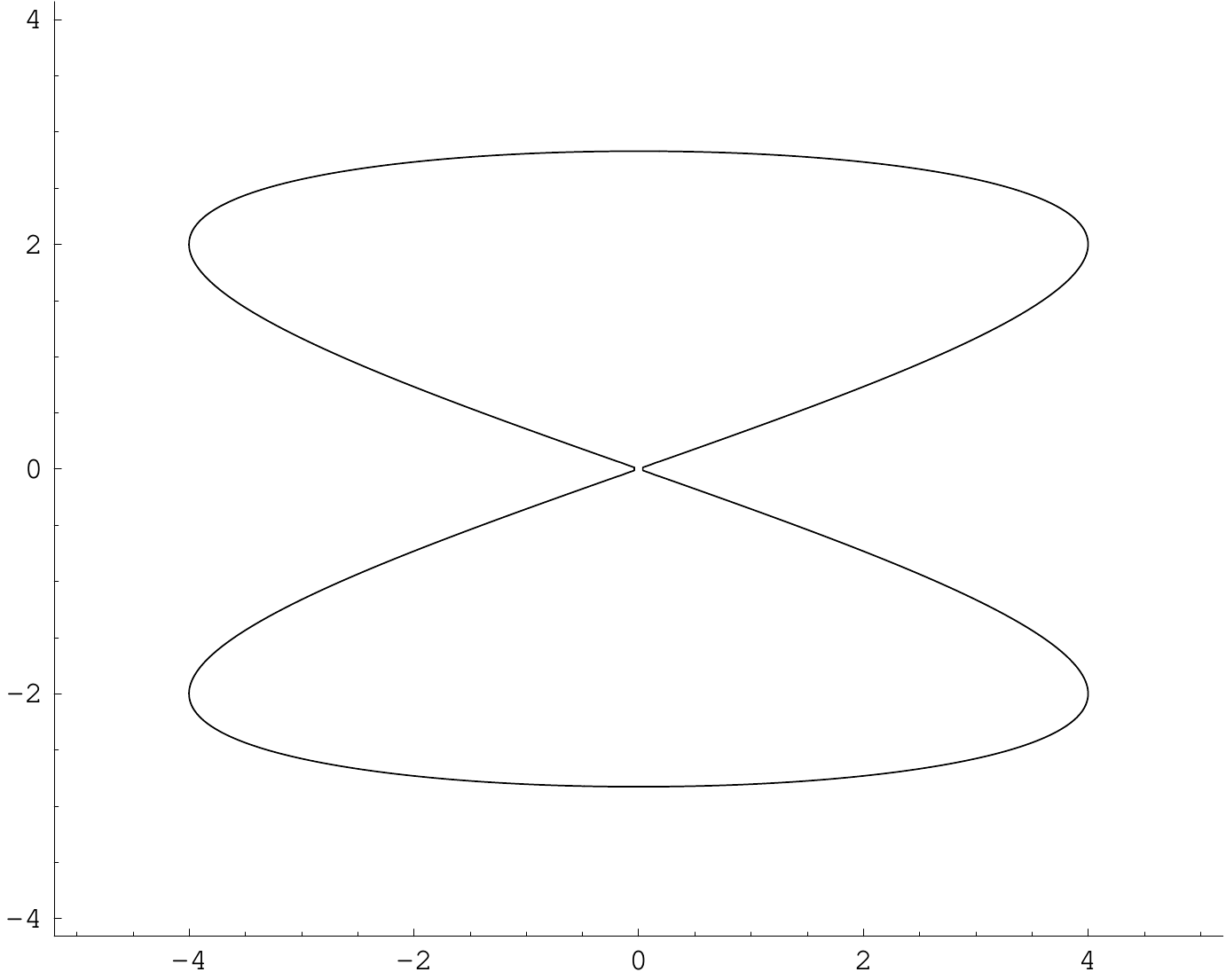}
\caption{Algebraic curve with denominator 2, the $1/2 \mapsto 1/4$ map.}
\label{mCurves_twos}
\end{center}
\end{figure}

\begin{figure}[h]
\begin{center}
   \includegraphics[width=2.5in]{Curves_1_4.pdf}     
   \includegraphics[width=2.5in]{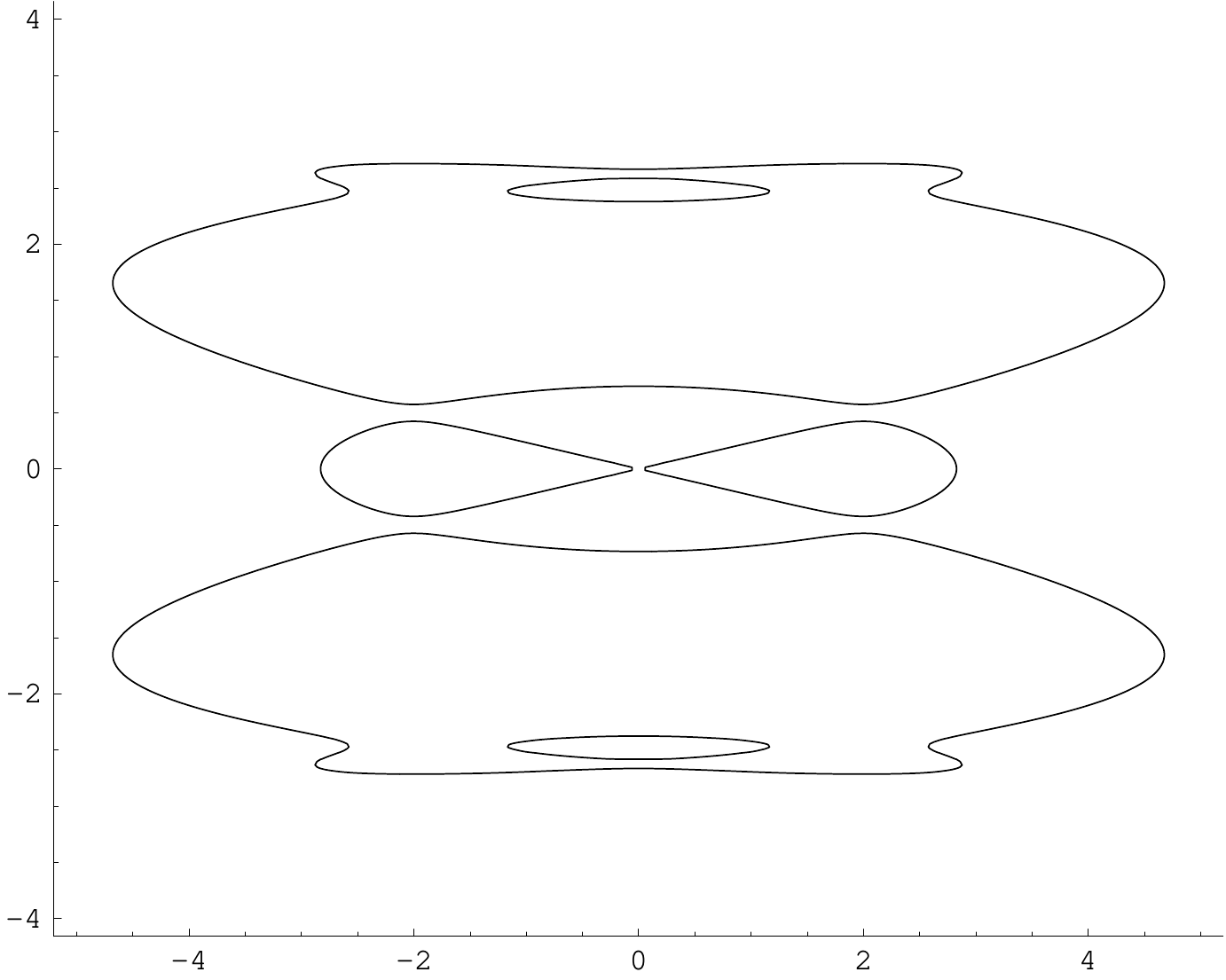} 
\caption{Algebraic curves with denominator 4: $1/4\mapsto 1/6$,  $3/4 \mapsto 3/10$. In each graph, there are four nearly linear segments down the diagonal of the graph. The right graph has some components which are disconnected from the linear segments.}
\label{mCurves_fours}
\end{center}
\end{figure}

.. \hfill ..
\pagebreak

\begin{figure}[h]
\begin{center}
$\begin{array}{ll }
   \includegraphics[width=2.5in]{Curves_1_6.pdf}   &  
   \includegraphics[width=2.5in]{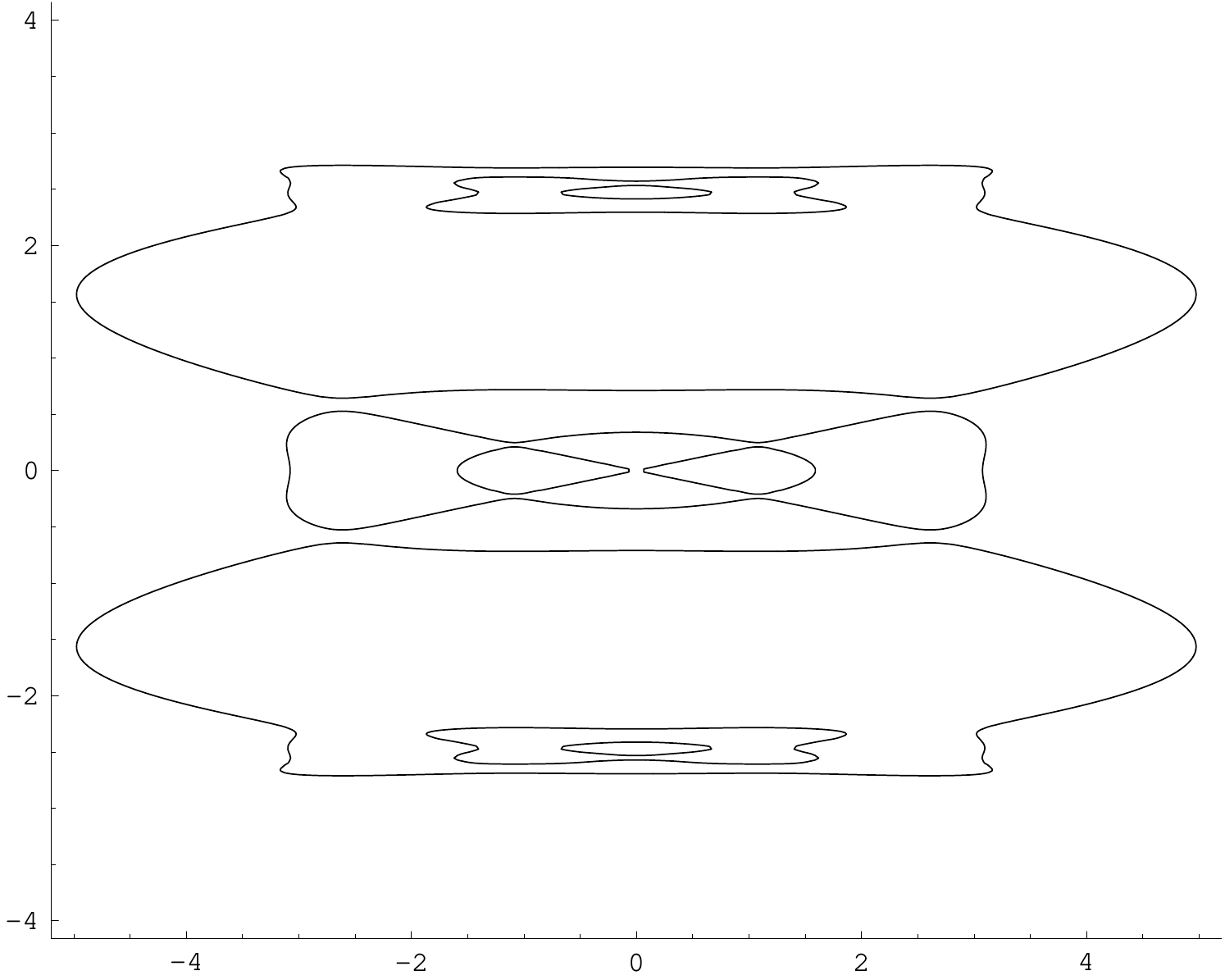}  \\
\end{array}
$
\caption{Algebraic curves with denominator 6: $1/6\mapsto 1/8$,  $5/6 \mapsto 5/16$. In each graph, there are six nearly linear segments down the diagonal of the graph. The right graph has some components which are disconnected from the linear segments.}
\label{mCurves_sixs}
\end{center}
\end{figure}

.. \hfill ..
\pagebreak

\begin{figure}[ht]
\begin{center}
$\begin{array}{ll }
   \includegraphics[width=2.5in]{Curves_1_8.pdf}   &  
   \includegraphics[width=2.5in]{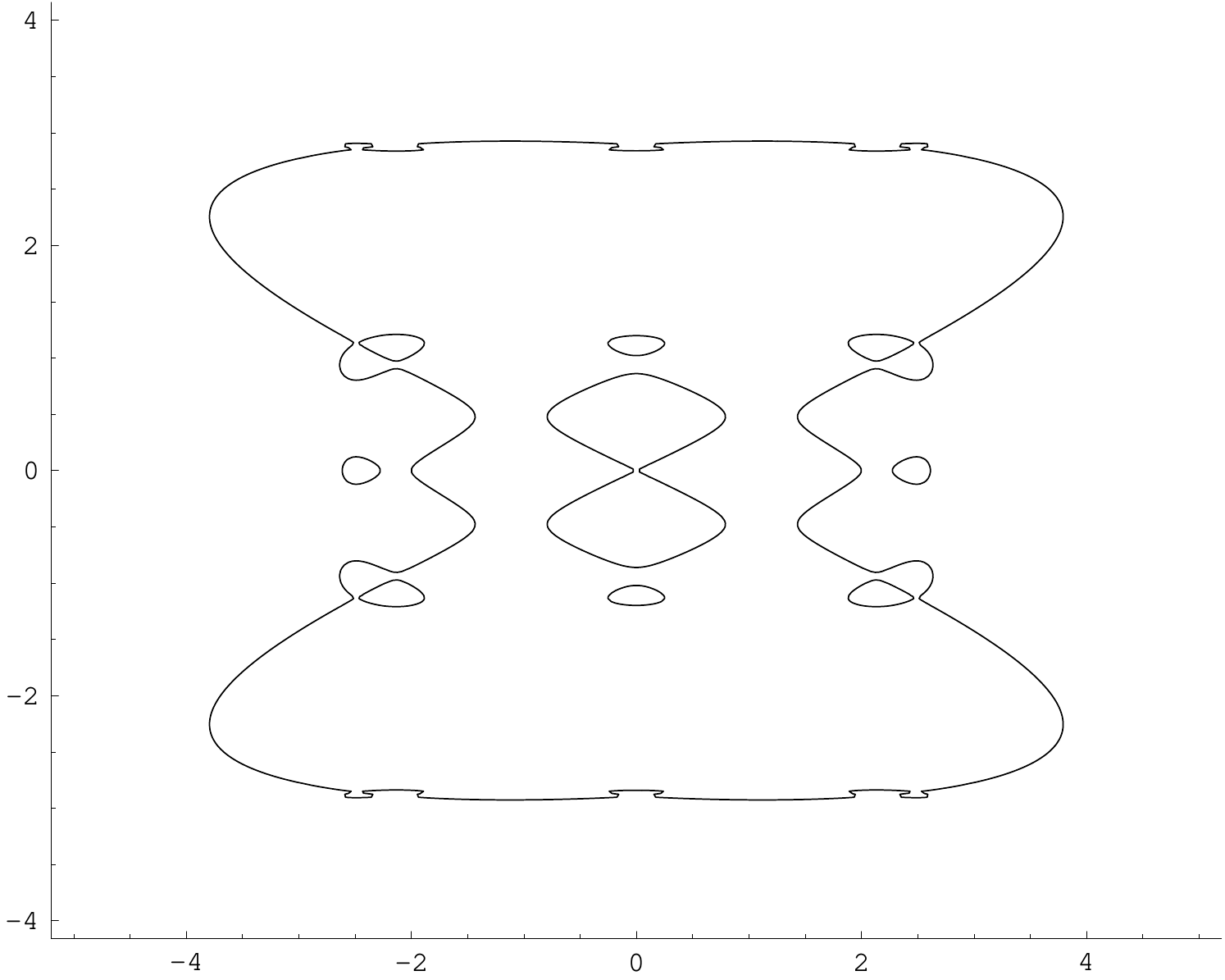}  \\
   \includegraphics[width=2.5in]{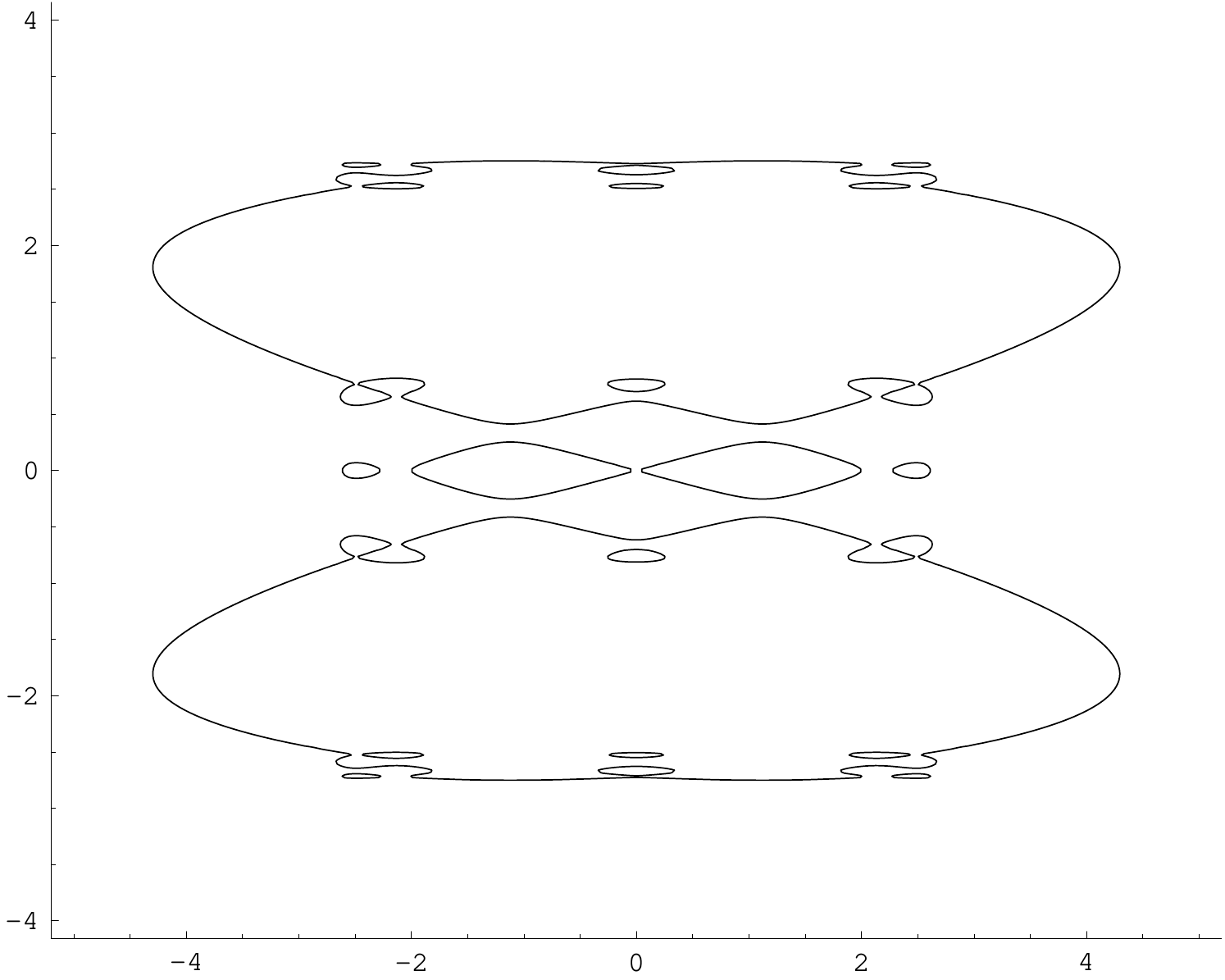}   &  
   \includegraphics[width=2.5in]{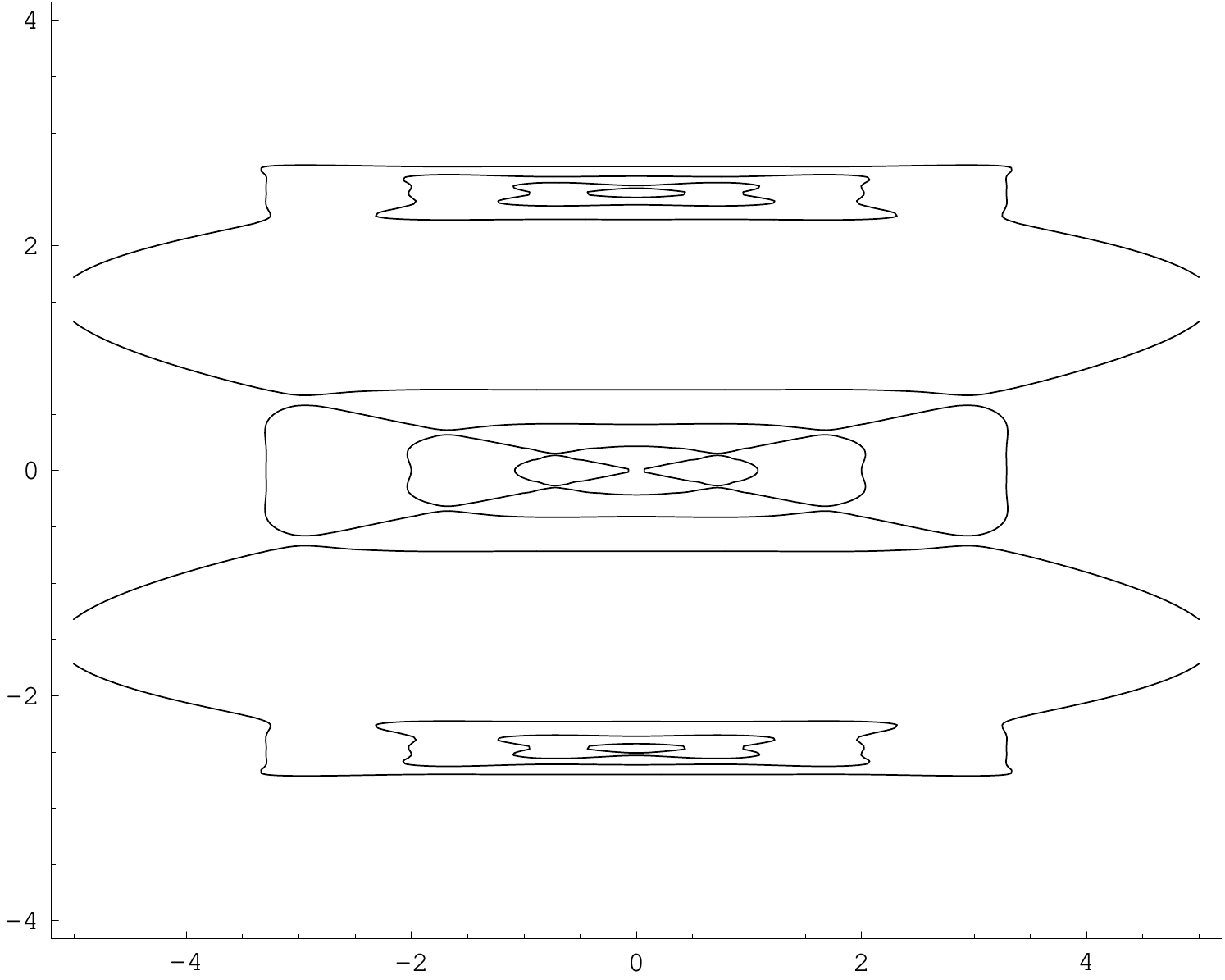}  \\
\end{array}
$
\caption{Algebraic curves with denominator 8: $1/8\mapsto 1/10$,  $3/8 \mapsto 3/14$,  $5/8 \mapsto 5/18$,  $7/8 \mapsto 7/22$. In each graph, there are eight nearly linear segments down the diagonal of the graph. All but the first graph has some components which are disconnected from the linear segments.}
\label{mCurves_eights}
\end{center}
\end{figure}

.. \hfill ..
\pagebreak

\end{document}